\documentclass[secthm]{elsart}

\journal{\ejc}

\usepackage[english]{babel}
\usepackage[T1]{fontenc}
\usepackage[latin1]{inputenc}

\usepackage{mathrsfs}
\usepackage{amsfonts}
\usepackage{amssymb}
\usepackage[fleqn]{amsmath}

\usepackage{comment}

\usepackage{pstricks}

\usepackage{color}
\usepackage{graphicx}

\usepackage[pagewise,mathlines]{lineno}

\DeclareMathAlphabet{\mathpzc}{OT1}{pzc}{m}{it}

\newcommand{\macro}[3]{\newcommand{#1}[#3]{#2}}

\newcommand\macros[4]
                   {%
                     \newenvironment{#1}[1][#2]{\par\vspace{1ex}
                       #3 \hspace{0.5em}#4}
                                    {\nopagebreak%
                                      \strut\nopagebreak%
                                      \par\vspace{2ex}
                                    }
                   }

\newtheorem{property}[thm]{\bfseries Property}

\macro{\twd}{\operatorname{twd}(#1)}{1}
\macro{\pwd}{\operatorname{pwd}(#1)}{1}
\macro{\bwd}{\operatorname{bwd}(#1)}{1}
\macro{\cwd}{\operatorname{cwd}(#1)}{1}
\macro{\mcwd}{\operatorname{mcwd}(#1)}{1}
\macro{\rwd}{\operatorname{rwd}(#1)}{1}
\macro{\sn}{\operatorname{sn}(#1)}{1}
\macro{\bn}{\operatorname{bn}(#1)}{1}

\macro{\brwd}{\operatorname{brwd}(#1)}{1}
\macro{\Frwd}{\operatorname{rwd}^{{#1}}(#2)}{2}
\macro{\Qrwd}{\Frwd{\gfq}{#1}}{1}
\macro{\frwd}{\Frwd{\bF}{#1}}{1}
\macro{\Fbrwd}{{#1}\textrm{-}\operatorname{brwd}(#2)}{2}
\macro{\fbrwd}{\Fbrwd{\bF}{#1}}{1}

\macro{\lrwd}{\operatorname{lrwd}(#1)}{1}
\macro{\lbwd}{\operatorname{lbwd}(#1)}{1}

\macro{\const}{\mathbf{#1}}{1}
\macro{\angl}{\mathop\langle #1 \mathop\rangle}{1}
\macro{\up}{\ulcorner #1\urcorner}{1}
\macro{\card}{\left|{#1}\right|}{1}
\macro{\floor}{\left\lfloor{#1}\right\rfloor}{1}
\macro{\ceil}{\left\lceil{#1}\right\rceil}{1}
\macro{\pare}{\left({#1}\right)}{1}
\macro{\crochet}{\left[{#1}\right]}{1}
\macro{\set}{\left\{{#1}\right\}}{1}
\macro{\range}{\set{{#1},\ldots,{#2}}}{2}
\macro{\mat}{M_{#1}}{1}
\macro{\matind}{{#1}[{#2},{#3}]}{3}
\macro{\matgind}{\matind{\matg}{#1}{#2}}{2}
\macro{\leaves}{\operatorname{L}_{#1}}{1}
\macro{\cutrk}{\operatorname{cutrk}^{{#1}}}{1}
\macro{\comp}{(X^{#1},#2\backslash X^{#1})}{2}
\macro{\supp}{\mathpzc{u}(#1)}{1}
\macro{\field}{\mathbb{F}_{#1}}{1}
\macro{\email}{\texttt{#1}}{1}
\macro{\subg}{#1\textrm{-}#2}{2}
\macro{\bicutrk}{\operatorname{bicutrk}^{{#1}}}{1}
\macro{\spivot}{\operatorname{\bar{*}}\limits_{#1}}{1}
\macro{\minor}{\operatorname{\parallel}\limits_{#1}}{1}

\def\gfa{\mathbb{a}}
\def\gfb{{\gfa^2}}
\def\gfq{\operatorname{\field{4}}}
\def\rk{\operatorname{rk}}

\def\matg{\mat{G}}
\def\fcutrk{\cutrk{\bF}}

\def\restriction#1#2{\mathchoice
              {\setbox1\hbox{${\displaystyle #1}_{\scriptstyle #2}$}
              \restrictionaux{#1}{#2}}
              {\setbox1\hbox{${\textstyle #1}_{\scriptstyle #2}$}
              \restrictionaux{#1}{#2}}
              {\setbox1\hbox{${\scriptstyle #1}_{\scriptscriptstyle #2}$}
              \restrictionaux{#1}{#2}}
              {\setbox1\hbox{${\scriptscriptstyle #1}_{\scriptscriptstyle #2}$}
              \restrictionaux{#1}{#2}}}
\def\restrictionaux#1#2{{#1\,\smash{\vrule height .8\ht1 depth .85\dp1}}_{\,#2}} 

\def\ie{\emph{i.e.}}

\def\bN{\mathbb{N}}

\def\bF{\mathbb{F}}

\def\bN{\mathbf{N}}

\def\bK{\mathbb{K}}

\def\cB{\mathcal{B}} 
\def\cS{\mathcal{S}} 

\def\cF{\mathcal{F}}

\def\cB{\mathcal{B}}

\def\cL{\mathcal{L}}

\def\cM{\mathcal{M}}

\def\tG{\widetilde{G}}

\def\dam{Discrete Applied Mathematics}

\def\jctb{Journal of Combinatorial Theory, Series B}

\def\jcss{Journal of Computer and System Sciences}

\def\siamjc{SIAM Journal on Computing}
\def\siamjdm{SIAM Journal on Discrete Mathematics}

\def\ejc{European Journal of Combinatorics}

\def\lap{Linear Algebra and its Applications}
\def\mp{Mathematical Programming}

\macros{defi}{Example}{\noindent{\bf #1}}{}
\macros{mthm}{Main Theorem}{\noindent{\bf #1}}{\itshape}

\newcounter{CaseCtr}
               {\end{list}}

\usepackage{calc}

{\end{list}}

\newcommand*\patchAmsMathEnvironmentForLineno[1]{%
  \expandafter\let\csname old#1\expandafter\endcsname\csname #1\endcsname
  \expandafter\let\csname oldend#1\expandafter\endcsname\csname end#1\endcsname
  \renewenvironment{#1}%
     {\linenomath\csname old#1\endcsname}%
     {\csname oldend#1\endcsname\endlinenomath}}%
\newcommand*\patchBothAmsMathEnvironmentsForLineno[1]{%
  \patchAmsMathEnvironmentForLineno{#1}%
  \patchAmsMathEnvironmentForLineno{#1*}}%
\AtBeginDocument{%
\patchBothAmsMathEnvironmentsForLineno{equation}%
\patchBothAmsMathEnvironmentsForLineno{align}%
\patchBothAmsMathEnvironmentsForLineno{flalign}%
\patchBothAmsMathEnvironmentsForLineno{alignat}%
\patchBothAmsMathEnvironmentsForLineno{gather}%
\patchBothAmsMathEnvironmentsForLineno{multline}%
}
\def\tsig{\tilde{\sigma}}

\begin{document}

\begin{frontmatter}

\title{Well-Quasi-Ordering of Matrices under Schur Complement and
  Applications to Directed Graphs}

\author{Mamadou Moustapha Kant\'e}

\address{Clermont-Université, Université Blaise Pascal, LIMOS,
  CNRS\\Complexe Scientifique des Cézeaux 63173 Aubiére Cedex,
  France\\ \email{mamadou.kante@isima.fr}}

\begin{abstract} 
  In [Rank-Width and Well-Quasi-Ordering of Skew-Symmetric or
    Symmetric Matrices, arXiv:1007.3807v1] Oum proved that, for a
  fixed finite field $\bF$, any infinite sequence $M_1,M_2,\ldots$ of
  (skew) symmetric matrices over $\bF$ of bounded \emph{$\bF$-rank-width} has a pair
  $i< j$, such that $M_i$ is isomorphic to a principal submatrix of a
  \emph{principal pivot transform} of $M_j$. We generalise this
  result to \emph{$\sigma$-symmetric matrices} introduced by Rao and
  myself in [The Rank-Width of Edge-Coloured Graphs,
    arXiv:0709.1433v4].  (Skew) symmetric matrices are special cases
  of $\sigma$-symmetric matrices. As a by-product, we obtain that for every infinite
  sequence $G_1,G_2,\ldots$ of directed graphs of bounded
  rank-width there exist a pair $i<j$ such that $G_i$ is a
  \emph{pivot-minor} of $G_j$. Another consequence is that non-singular
  principal submatrices of a $\sigma$-symmetric matrix form a
  \emph{delta-matroid}. We extend in this way the notion of
  representability of delta-matroids by Bouchet.
\end{abstract}

\begin{keyword} rank-width; sigma-symmetry; 
  edge-coloured graph; well-quasi-ordering; principal pivot transform;
  pivot-minor.
\end{keyword}

\end{frontmatter}

\section{Introduction}\label{sec:1}



\emph{Clique-width} \cite{CER93} is a graph complexity measure that
emerges in the works by Courcelle et al. (see for instance the book
\cite{CE11}). It extends \emph{tree-width} \cite{RS90} in the sense
that graph classes of bounded tree-width have bounded clique-width,
but the converse is false (distance hereditary graphs have
clique-width at most $3$ and unbounded tree-width). Clique-width has
similar algorithmic properties as tree-width and seems to be the right
complexity measure for the investigations of polynomial time
algorithms in dense graphs for a large set of NP-complete problems
\cite{CE11}. It is then important to identify graph classes of bounded
clique-width. Unfortunately, contrary to tree-width, there is no known
polynomial time algorithm that checks if a given graph has
clique-width at most $k$, for fixed $k\geq 4$ (for $k\leq 3$, see the
algorithm by Corneil et al.  \cite{CHLRR00}). Furthermore,
clique-width is not monotone with respect to \emph{graph minor} (cliques have
clique-width $2$) and is only known to be monotone with respect to the
\emph{induced subgraph} relation which is not a well-quasi-order on
graph classes of bounded clique-width (cycles have clique-width at
most $4$ and are not well-quasi-ordered by the induced subgraph
relation).

In their investigations for a recognition algorithm for graphs of
clique-width at most $k$, for fixed $k$, Oum and Seymour \cite{OS06}
introduced the complexity measure \emph{rank-width} of undirected
graphs. Rank-width and clique-width of undirected graphs are
equivalent in the sense that a class of undirected graphs has bounded
rank-width if and only if it has bounded clique-width. But, if
rank-width shares with clique-width its same algorithmic properties
(see for instance \cite{CK09}), it has better combinatorial
properties. 
\begin{enumerate}
\item There exists a cubic-time algorithm that checks whether an
  undirected graph has rank-width at most $k$, for fixed $k$
  \cite{HO07}. 

\item Rank-width is monotone with respect to the \emph{pivot-minor}
  relation. This relation generalises the notion of graph minor
  because if $H$ is a minor of $G$, then $I(H)$, the \emph{incidence
    graph} of $H$, is a pivot-minor of $I(G)$. Undirected graphs of
  rank-width at most $k$ are characterised by a finite list of
  undirected graphs to exclude as pivot-minors \cite{OUM05}.

\item Furthermore, rank-width is related to the \emph{branch-width} of
  binary matroids. Branch-width of matroids plays an important role in
  the project by Geelen et al. \cite{GGW07b} aiming at extending
  techniques in the Graph Minors Project to representable matroids
  over finite fields in order to prove that representable matroids
  over finite fields are well-quasi-ordered by \emph{matroid
    minors}. Such a result would answer positively Rota's Conjecture
  \cite{GGW07b}. It turns out that the branch-width of a binary
  matroid is one more than the rank-width of its fundamental graphs
  and a fundamental graph of a minor of a matroid $\cM$ is a
  pivot-minor of a fundamental graph of $\cM$.
\end{enumerate}

It is then relevant to ask whether undirected graphs are
well-quasi-ordered by the pivot-minor relation. This would imply that
binary matroids are well-quasi-ordered by matroid minors, and hence
the \emph{Graph Minor Theorem} \cite{RS04}. This would also help understand the structure
of graph classes of bounded clique-width and of many dense graph
classes where the Graph Minor Theorem  fails to explain
their structure. Geelen et al. have successfully adapted many
techniques in the Graph Minors Project \cite{RS} and obtained
generalisations of some results in the Graph Minors Projects to
representable matroids over finite fields (see the survey
\cite{GGW07b}). Inspired by the links between rank-width and
branch-width of binary matroids, Oum \cite{OUM08} adapted the
techniques by Geelen et al. and proved that undirected graphs of
bounded rank-width are well-quasi-ordered by the pivot-minor relation.
As for the Graph Minors Project, this seems to be a first step towards
a Graph Pivot-Minor Theorem.

However, rank-width has a drawback: it is defined in Oum's works only
for undirected graphs. But, clique-width was originally defined for
graphs (directed or not, with edge-colours or not). Hence, one would
know about the structure of (edge-coloured) directed graphs of bounded
clique-width. Rao and myself \cite{KR11} we have defined a notion of
rank-width, called \emph{$\bF$-rank-width}, for $\bF^*$-graphs, \ie,
graphs with edge-colours from a field $\bF$, and explained how to use
it to define a notion of rank-width for graphs (directed or not, with
edge-colours or not). Moreover, the notion of rank-width of undirected
graphs is a special case of it.  $\bF$-rank-width is equivalent to
clique-width and all the known results, but the well-quasi-ordering
theorem by Oum \cite{OUM08}, concerning the rank-width of undirected
graphs have been generalised to the $\bF$-rank-width of
$\bF^*$-graphs. We complete the tableau in this paper by proving a
well-quasi-ordering theorem for $\bF^*$-graphs of bounded
$\bF$-rank-width, and hence for directed graphs.

In  \cite{OUM10} Oum noticed that the \emph{principal pivot transform}
introduced by Tucker \cite{TUC60} can be used to obtain a
well-quasi-ordering theorem for (skew) symmetric matrices over finite
fields of bounded $\bF$-rank-width.  This result unifies his own
result on the well-quasi-ordering of undirected graphs of bounded
rank-width by pivot-minor\cite{OUM08}, the well-quasi-ordering by
matroid minor of matroids representable over finite fields of bounded
branch-width \cite{GGW02} and the well-quasi-ordering by graph minor
of undirected graphs of bounded tree-width \cite{RS90}. In order to
prove the well-quasi-ordering theorem for $\bF^*$-graphs of bounded
$\bF$-rank-width, we will adapt the techniques used by Oum in
\cite{OUM10} to \emph{$\sigma$-symmetric} matrices. The notion of
$\sigma$-symmetric matrices were introduced by Rao and myself in
\cite{KR11} and subsumes the notion of (skew) symmetric
matrices. Oum's proof can be summarised into two steps.
\begin{itemize}
\item[(i)] He first developed a theory about the notion of
  \emph{lagrangian chain-groups}, which are generalisations of
  \emph{isotropic systems} \cite{BOU87a} and of \emph{Tutte
    chain-groups} \cite{TUT71}. Tutte chain-groups are another
  characterisation of representable matroids, and isotropic systems
  are structures that extend some properties of $4$-regular graphs and
  of circle graphs. Isotropic systems played an important role in the
  proof of the well-quasi-ordering of undirected graphs of bounded
  rank-width by pivot-minor. As for Tutte chain groups and isotropic
  systems, lagrangian chain groups are vector spaces equipped with a
  bilinear form. Oum introduced a notion of minor for lagrangian chain
  groups that subsumes the matroid minor and the notion of minor of
  isotropic systems. He also defined a connectivity function for
  lagrangian chain groups that generalises the connectivity function
  of matroids and allows to define a notion of \emph{branch-width} for
  them. He then proved that lagrangian chain-groups of bounded
  branch-width are well-quasi-ordered by lagrangian chain groups
  minor.

\item[(ii)] He secondly proved that to any lagrangian chain-group, one
  can associate a (skew) symmetric matrix and vice-versa. These
  matrices are called \emph{matrix representations} of lagrangian
  chain-groups. He can thus formulate the well-quasi-ordering theorem
  of lagrangian chain-groups in terms of (skew) symmetric matrices.
\end{itemize}

We will follow the same steps. We will extend the notion of lagrangian
chain-groups to make it compatible with $\sigma$-symmetric
matrices. Then, we prove that these lagrangian chain-groups admit
representations by $\sigma$-symmetric matrices. 

The paper is organised as follows. We present some notations needed
throughout the paper in Section \ref{sec:2}. Chain groups are
revisited in Section \ref{sec:3}. Section \ref{sec:4} is devoted to
the links between chain groups and $\sigma$-symmetric matrices. The
main theorem (Theorem \ref{thm:4.2}) of the paper is presented in
Section \ref{sec:4}. Applications to directed graphs and more
generally to edge-coloured graphs is presented in Section
\ref{sec:5}. An old result by Bouchet \cite{BOU88} states that
non-singular principal submatrices of a (skew) symmetric matrix form
a \emph{delta-matroid}. We extend this result to $\sigma$-symmetric
matrices and obtain a new notion of representability of
delta-matroids in Section \ref{sec:6}.   

\section{Preliminaries}\label{sec:2}

For two sets $A$ and $B$, we let $A\setminus B$ be the set $\{x\in
A\mid x\notin B\}$. The power-set of a set $V$ is denoted by $2^V$. We
often write $x$ to denote the set $\{x\}$. We denote by $\bN$ the set
containing zero and the positive integers. If $f:A\to B$ is a
function, we let $\restriction{f}{X}$, the restriction of $f$ to
$X\subseteq A$, be the function $\restriction{f}{ X}:X\to B$ where for
every $a\in X,\ \restriction{f}{ X}(a) := f(a)$. For a finite set $V$,
we say that the function $f:2^V\to \bN$ is \emph{symmetric} if for any
$X\subseteq V, ~f(X)=f(V\setminus X)$; $f$ is \emph{submodular} if
for any $X,Y\subseteq V$, $f(X\cup Y) + f(X\cap Y) \leq f(X) +f(Y)$.

We denote by $+$ and $\cdot$ the binary operations of any field and by
$0$ and $1$ the identity elements of $+$ and $\cdot$ respectively.
Fields are denoted by the symbol $\bF$ and finite fields of order $q$
by $\bF_q$. We recall that finite fields are commutative. For a field
$\bF$, we let $\bF^*$ be the set $\bF\setminus \{0\}$. We refer to
\cite{LN97} for our field terminology.

We use the standard graph terminology, see for instance
\cite{DIE05}. A \emph{directed graph} $G$ is a couple $(V_G,E_G)$
where $V_G$ is the set of vertices and $E_G\subseteq V_G\times V_G$ is
the set of edges. A directed graph $G$ is said to be \emph{undirected}
if $(x,y)\in E_G$ implies $(y,x)\in E_G$.  For a directed graph $G$,
we denote by $G[X]$, called the subgraph of $G$ induced by $X\subseteq
V_G$, the directed graph $(X,E_G\cap (X\times X))$. The degree of a
vertex $x$ in an undirected graph $G$ is the cardinal of the set
$\{y\mid xy\in E_G\}$.  Two directed graphs $G$ and $H$ are
\emph{isomorphic} if there exists a bijection $h:V_G\to V_H$ such that
$(x,y)\in E_G$ if and only if $(h(x),h(y))\in E_H$. We call $h$ an
\emph{isomorphism} between $G$ and $H$. All directed graphs are finite
and can have loops.

A \emph{tree} is an acyclic connected undirected graph. A \emph{cubic
  tree} is a tree such that the degree of each vertex is either $1$ or
$3$. For a tree $T$ and an edge $e$ of $T$, we let $\subg{T}{e}$
denote the graph $(V_T,E_T\setminus \{e\})$.

A \emph{layout} of a finite set $V$ is a pair $(T,\cL)$ of a cubic
tree $T$ and a bijective function $\cL$ from the set $V$ to the set
$\leaves{T}$ of vertices of degree $1$ in $T$. For each edge $e$ of
$T$, the connected components of $\subg{T}{e}$ induce a bipartition
$(X_e,V\setminus X_e)$ of $\leaves{T}$, and thus a bipartition
$\comp{e}{V}= (\cL^{-1}(X_e), \cL^{-1}(V\setminus X_e))$ of $V$. Let
$f:2^V\to \bN$ be a symmetric function and $(T,\cL)$ a layout of
$V$. The \emph{$f$-width of each edge $e$ of $T$} is defined as
$f(X^e)$ and the \emph{$f$-width of $(T,\cL)$} is the maximum
$f$-width over all edges of $T$. The \emph{$f$-width of $V$} is the
minimum $f$-width over all layouts of $V$. The notions of layout and
of $f$-width are commonly called \emph{branch-decomposition} and
\emph{branch-width} of $f$. However, this terminology is not
appropriate since $f$ is only a measure for the cuts $(\cL^{-1}(X_e),
\cL^{-1}(V\setminus X_e))$ and other measures could be used with the
same layout.

\subsection{Well-Quasi-Order}\label{subsec:2.3}
We review in this section the \emph{well-quasi-ordering} notion. A
binary relation is a \emph{quasi-order} if it is reflexive and
transitive. A quasi-order $\preceq$ on a set $\mathcal{U}$ is a
\emph{well-quasi-order}, and the elements of $\mathcal{U}$ are
\emph{well-quasi-ordered} by $\preceq$, if for every infinite sequence
$x_0,x_1,\ldots$ in $\mathcal{U}$ there exist $i<j$ such that $x_i\preceq
x_j$. The notion of well-quasi-ordering is flourishing and there exist
several equivalent definitions of the well-quasi-ordering notion. For
instance, a quasi-order $\preceq$ on a set $\mathcal{U}$ is a
well-quasi-order if and only if $\mathcal{U}$ contains no infinite
antichain and no infinite strictly decreasing sequence. One
consequence of this characterisation is that every $\preceq$-closed set
$X$ of $\mathcal{U}$, \ie, if $y\in X$ and $x\preceq y$ then $x\in X$, is
characterised by a finite list $Forb(X)$ such that $x\in X$ if and
only if there is no $z\in Forb(X)$ with $z\preceq x$. Hence, the
well-quasi-ordering notion is an interesting tool for characterising
graph classes. There exist several well-quasi-ordering theorems in
the literature, see for instance \cite[Chapter 12]{DIE05} for some of
them.

\subsection{Sesqui-Morphism}\label{subsec:2.1}

We recall the notion of \emph{sesqui-morphism} introduced in
\cite{KR11} in order to extend the notion of rank-width to directed
graphs. Let $\bF$ be a field and $\sigma:\bF\to \bF$ a bijection. We
recall that $\sigma$ is an involution if $\sigma\circ \sigma$ is the
identity. We call $\sigma$ a \emph{sesqui-morphism} if $\sigma$ is an
involution, and the function $\tsig:=[x\mapsto \sigma(x)/\sigma(1)]$
is an automorphism.  It is worth noticing that if $\sigma:\bF\to \bF$
is a sesqui-morphism, then $\sigma(0)=0$ and for every $a,b\in \bF$,
$\sigma(a+b)=\sigma(a)+\sigma(b)$. Moreover, $\tsig$ is an
involution. The next proposition summarises some properties of
sesqui-morphisms.

\begin{prop}\label{prop:2.1} Let $\sigma:\bF\to \bF$ be a
  sesqui-morphism. Then, for all $a,b,a_i\in \bF$, $c\in \bF^*$ and
  all $n\in \bN$,
  \begin{align}
\sigma(-a) &= -\sigma(a) \\
\sigma(a_1\cdot a_2 \cdots  a_n) & = \frac{\sigma(a_1)\cdot
  \sigma(a_2) \cdots  \sigma(a_n)}{\sigma(1)^{n-1}} &
\\
\sigma(a^n) & = \frac{\sigma(a)^n}{\sigma(1)^{n-1}} \\
\sigma(a^{-n}) &= \frac{\sigma(1)^{n+1}}{\sigma(a)^n}\\
\sigma\left(\frac{a}{c}\right) &= \frac{\sigma(1)\cdot
  \sigma(a)}{\sigma(c)}\\
\sigma\left(\frac{a\cdot b}{c}\right) & = \frac{\sigma(a)\cdot
  \sigma(b)}{\sigma(c)} 
  \end{align}
\end{prop}

\begin{pf*}{Proof.} Equation (1) is trivial since
  $\sigma(a)+\sigma(-a)=\sigma(a-a)= \sigma(0)=0$.

  Equation (2) will be proved by induction. The case $n=2$ is trivial
  since $\tsig$ is an automorphism. Assume $n>2$. Then, 
  \begin{align*}
    \sigma(a_1\cdot a_2\cdots a_n) & = \sigma(a_1\cdot a_2 \cdots
    a_{n-1}) \cdot \frac{\sigma(a_n)}{\sigma(1)}\\ &= \frac{\sigma(a_1) \cdot
    \sigma(a_2) \cdots \sigma(a_{n-1})}{\sigma(1)^{n-2}} \cdot
    \frac{\sigma(a_n)}{\sigma(1)} 
  \end{align*}
  This proves the equation. Equation (3) is a direct consequence of
  Equation (2) since $\sigma(a^n) = \sigma(\underbrace{a\cdots a}_{n})$. 

  Since $\sigma(a^{-n}) = \tsig(a^{-n}) \cdot \sigma(1)$, Equation (4)
  follows from this equality $\tsig(a^{-n})= \frac{1}{\tsig(a^n)}$. Equations
  (5) and (6) are consequences of Equations (2)-(4). \qed
\end{pf*}

Examples of sesqui-morphisms are the identity automorphism (called
\emph{symmetric sesqui-morphism}) and the function $[x\mapsto -x]$
(called \emph{skew-symmetric sesqui-morphism}). The next proposition
states that they are the only ones in prime fields.

\begin{prop}\label{prop:2.2} Let $p$ be a prime number and let
  $\sigma:\bF_p\to \bF_p$ be a function. Then, $\sigma$ is a
  sesqui-morphism if and only if $\sigma$ is symmetric or
  skew-symmetric.
\end{prop}

\begin{pf*}{Proof.} Assume $\sigma:\bF_p\to \bF_p$ is a
  sesqui-morphism. It is well-known that the only automorphism in
  $\bF_p$, $p$ prime, is the identity. Hence, $\tsig(a)=a$ for all
  $a\in \bF_p$. Thus, $\sigma(a)=a\cdot \sigma(1)$, and hence, $1=
  \sigma(\sigma(1)) = \sigma(1)^2$. Therefore, $\sigma(1)=\pm 1$.\qed
\end{pf*}

Along this paper, sesqui-morphisms will be denoted by the Greek letter
$\sigma$, and then we will often omit to say "let $\sigma:\bF\to \bF$
be a sesqui-morphism". 

\subsection{Matrices and $\bF$-Rank-Width} \label{subsec:2.2}


For sets $R$ and $C$, an \emph{$(R,C)$-matrix} is a matrix where the
rows are indexed by elements in $R$ and columns indexed by elements in
$C$. If the entries are over a field $\bF$, we call it an
$(R,C)$-matrix over $\bF$. For an $(R,C)$-matrix $M$, if $X\subseteq
R$ and $Y\subseteq C$, we let $\matind{M}{X}{Y}$ be the submatrix of
$M$ where the rows and the columns are indexed by $X$ and $Y$
respectively. Along this paper matrices are denoted by capital
letters, which will allow us to write $m_{xy}$ for $\matind{M}{x}{y}$
when it is possible. The matrix rank-function is denoted $\rk$.  We
will write $M[X]$ instead of $M[X,X]$ and such submatrices are called
\emph{principal submatrices}.  The transpose of a matrix $M$ is
denoted by $M^t$, and the inverse of $M$, if it exists, \ie, if $M$ is
\emph{non-singular}, is denoted by $M^{-1}$. The \emph{determinant} of
$M$ is denoted by $\det(M)$.  A $(V_1,V_1)$-matrix $M$ is said
\emph{isomorphic} to a $(V_2,V_2)$-matrix $N$ if there exists a
bijection $h:V_1\to V_2$ such that $m_{xy}=n_{\scriptsize
  h(x)h(y)}$. We refer to \cite{LIP91} for our linear algebra
terminology.

For a sesqui-morphism $\sigma:\bF\to \bF$, a $(V,V)$-matrix $M$ over
$\bF$ is said \emph{$\sigma$-symmetric} if $m_{yx} = \sigma(m_{xy})$
for all $x,y\in V$.  Examples of $\sigma$-symmetric matrices are
(skew) symmetric matrices with $\sigma$ being the (skew) symmetric
sesqui-morphism. From Proposition \ref{prop:2.2} they are the only
$\sigma$-symmetric matrices over prime fields. A $(V,V)$-matrix $M$ is
said \emph{$(\sigma,\epsilon)$-symmetric} if $\epsilon(x)\cdot m_{xy}
= \epsilon(y)\cdot \sigma(m_{yx})$ for all $x,y\in V$, $\epsilon:V\to
\{-1,+1\}$ being a function. If $\sigma$ is the (skew) symmetric
sesqui-morphism, $(\sigma,\epsilon)$-matrices are called matrices of
\emph{symmetric type} in \cite{BOU88}. It is worth noticing that a
matrix is $\sigma$-symmetric if and only if it is
$(\sigma,\epsilon)$-symmetric with $\epsilon$ a constant function.

We recall now the notion of \emph{$\bF$-rank-width} of
$(\sigma,\epsilon)$-matrices. It will be used to extend the notion of
rank-width to directed graphs. The \emph{$\bF$-cut-rank} function of a
$(\sigma,\epsilon)$-symmetric $(V,V)$-matrix $M$ is the function
$\fcutrk_M:2^{V}\to \bN$ where $\fcutrk_M(X) =
\rk(\matind{M}{X}{V\setminus X})$ for all $X\subseteq V$. From
Proposition \ref{prop:3.5} and Theorem \ref{thm:4.1}, the function
$\fcutrk_M$ is symmetric and submodular (a more direct proof for
$\sigma$-symmetric matrices can be found in \cite{KR11}, but it can be
easily adapted to $(\sigma,\epsilon)$-symmetric matrices). The
\emph{$\bF$-rank-width} of a $(\sigma,\epsilon)$-symmetric
$(V,V)$-matrix $M$ is the $\fcutrk_M$-width of $V$.

If $G$ is an undirected graph, then its adjacency matrix $A_G$ over
$\bF_2$ is $\sigma_1$-symmetric, with $\sigma_1$ the identity
automorphism on $\bF_2$. One easily checks that the rank-width of $G$
\cite{OUM05} is exactly the $\bF_2$-rank-width of $A_G$. 

Let $M$ be a matrix of the form $\left(\begin{smallmatrix} A & B \\ C
  & D \end{smallmatrix}\right)$ where $A:=M[X]$ is non-singular. The
\emph{Schur complement of $A$ in $M$}, denoted by $M/A$, is $D - C
\cdot A^{-1} \cdot B$. Oum proved the following.

\begin{thm}[\cite{OUM10}] \label{thm:2.1} Let $\bF$ be a finite field and $k$ a
  positive integer. For every infinite sequence $M_1,M_2,\ldots$ of
  symmetric or skew-symmetric matrices over $\bF$ of $\bF$-rank-width
  at most $k$, there exist $i<j$ such that $M_i$ is isomorphic to a
  principal submatrix of $M_j/A$ for some non-singular principal
  submatrix $A$ of $M_j$.
\end{thm}

This theorem unifies in a single one the well-quasi-ordering theorems
in \cite{GGW02,OUM08,RS90}. We will show that this theorem still holds
in the case of $(\sigma,\epsilon)$-symmetric matrices that are not
necessarily (skew) symmetric. As a by product, we will get a
well-quasi-ordering theorem for directed graphs. In order to do so, we
will adapt the same techniques as Oum's proof.

\section{Chain Groups Revisited} \label{sec:3}


\emph{Chain groups} were introduced by Tutte \cite{TUT71} for matroids
and were also studied by Bouchet in his series of papers dealing with
circle graphs and eulerian circuits of $4$-regular graphs (see for
instance \cite{BOU87a,BOU87b,BOU88}). 

The key point in the proof of Theorem \ref{thm:2.1} is to associate to
each (skew) symmetric matrix a chain group and then use the
well-quasi-ordering theorem on chain groups. We will revise the
definitions by Oum so that to associate to each
$(\sigma,\epsilon)$-symmetric matrix a chain group. All the vector
spaces manipulated have finite dimension.  The dimension of a vector
space $W$ is denoted by $\dim(W)$. If $f:W\to V$ is a linear
transformation, we denote by $Ker(f)$ the set $\{u\in W\mid f(u)=0\}$
and $Im(f)$ the set $\{f(u)\in V\mid u\in W\}$. It is worth noticing
that both are vector spaces. For a vector space $K$,
we let $K^*:=K\setminus \{0\}$.

For a field $\bF$ and sesqui-morphism $\sigma:\bF\to \bF$, we let
$\bK_\sigma$ be the $2$-dimensional vector space $\bF^2$ over $\bF$
equipped with the application $\mathbf{b}_\sigma:\bK_\sigma\times \bK_\sigma \to
\bF$ where $\mathbf{b}_\sigma( \left( \begin{smallmatrix}
  a\\ b \end{smallmatrix} \right), \left( \begin{smallmatrix}
  c\\ d \end{smallmatrix} \right) ) = \sigma(1)\cdot a \cdot \sigma(d)
-b \cdot \sigma(c)$.  The application $\mathbf{b}_\sigma$ is not
bilinear, however it is linear with respect to its left operand, which
is enough for our purposes.  It is worth noticing that if $\sigma$ is
skew-symmetric (or symmetric), then $\mathbf{b}_\sigma$ is what is
called $b^+$ (or $b^-$) in \cite{OUM10}. The following properties are
easy to obtain from the definition of $\mathbf{b}_\sigma$.

\begin{property}\label{pt:3.1} Let $u,v, w\in \bK_\sigma$ and $k \in
  \bF$. Then, 
\begin{align*}
\mathbf{b}_\sigma(u+v,w) &= \mathbf{b}_\sigma(u,w) + \mathbf{b}_\sigma(v,w),\\
\mathbf{b}_\sigma(u, v+w) &= \mathbf{b}_\sigma(u,v) + \mathbf{b}_\sigma(u,w),\\
\mathbf{b}_\sigma(k\cdot u, v) &= k \cdot \mathbf{b}_\sigma(u,v),\\
\mathbf{b}_\sigma(u,k \cdot v) &= \tsig (k) \cdot
\mathbf{b}_\sigma(u,v).\\
\sigma(\mathbf{b}_\sigma(u,v)) &= \frac{-1}{\sigma(1)^2} \cdot \mathbf{b}_\sigma(v,u).
\end{align*}
\end{property}

\begin{property}\label{pt:3.2} Let $u\in \bK_\sigma$. 
  \begin{enumerate}
  \item[(i)] If $\mathbf{b}_\sigma(u,v) = 0$
  for all $v\in \bK_\sigma$, then $u=0$.

\item[(ii)] If $\mathbf{b}_\sigma(v,u)=0$ for all $v\in
  \bK_\sigma$, then $u=0$.
\end{enumerate}
\end{property}

Let $W$ be a vector space over $\bF$ and $\varphi:W\times W\to \bF$ a
function. If $\varphi$ satisfies equalities in Property \ref{pt:3.1},
we call it a \emph{$\sigma$-sesqui-bilinear form}.  It is called a
\emph{non-degenerate} $\sigma$-sesqui-bilinear form if it also
satisfies Property \ref{pt:3.2}.

Let $W$ be a vector space over $\bF$ equipped with $\varphi$ a
non-degenerate $\sigma$-sesqui-bilinear form. A vector $u$ is said
\emph{isotropic} if $\varphi(u,u)=0$. A subspace $L$ of $W$ is called
\emph{totally isotropic} if $\varphi(u,v)=0$ for all $u,v\in L$. For a
subspace $L$ of $W$, we let $L^{\bot}:=\{v\in W\mid \varphi(u,v)=0$
for all $u\in L\}$. It is worth noticing that if $L$ is totally
isotropic, then $L\subseteq L^\bot$. The following theorem is a
well-known theorem in the case where $\varphi$ is a non-degenerate
bilinear form.

\begin{thm}\label{thm:3.1} Let $W$ be a vector space over $\bF$
  equipped with a non-degenerate $\sigma$-sesqui-bilinear form $\varphi$. Then,
  $\dim(L) + \dim(L^\bot) = \dim(W)$ for any subspace $L$ of $W$.
\end{thm}

\begin{pf*}{Proof.} The proof is a standard one. We denote by $W^*$
  the set of linear transformations $[W\to \bF]$. It is well-known that
  $W^*$ is a vector space. Let $\varphi_R:W\to
  W^*$ such that $\varphi_R(u) := [w\mapsto \varphi(w,u)]$. From
  Property \ref{pt:3.1}, $\varphi_R$ is clearly a linear
  transformation. Let $\alpha$ be a restriction of $\varphi_R$ to $L$. By
  a well-known theorem in linear algebra, $\dim(L)= \dim(Ker(\alpha)) +
  \dim(Im(\alpha))$. 

  By definition, $Ker(\alpha) = \{u\in L\mid \varphi(w,u)=0$ for all
  $w\in W\}$, which is equal to $\{0\}$ since $\varphi$ is
  non-degenerate. Hence, $\dim(Ker(\alpha))=0$, \ie,
  $\dim(L)=\dim(Im(\alpha))$. 

  If we let $Im(\alpha)^{\circ}:=\{v\in W\mid \theta(v)=0$ for all $\theta
  \in Im(\alpha)\}$, we know by a theorem in linear algebra that
  $\dim(Im(\alpha))+\dim(Im(\alpha)^{\circ}) = \dim(W^*)$. But, 
  \begin{align*}
    Im(\alpha)^{\circ} &= \{v\in W\mid \alpha(w)(v) = 0\ \textrm{for
      all}\ w\in L\} \\ & = \{v\in W\mid \varphi(v,w) = 0\ \textrm{for
      all}\ w\in L\} = L^\bot.
  \end{align*}
  Hence, $\dim(L) = \dim(W^*) - \dim(L^\bot) = \dim (W) - \dim(L^\bot)$
  since $\dim(W^*) = \dim(W)$. \qed
\end{pf*}

As a consequence, we get that $L= (L^\bot)^\bot$. And, if $L$ is
totally isotropic, then $2\cdot \dim(L)\leq \dim(W)$.  

Let $V$ be a finite set and $K$ a vector space over $\bF$. A
\emph{$K$-chain on $V$} is a function $f:V\to K$. We let $K^V$ be the set
of $K$-chains on $V$. It is well-known that $K^V$ is a vector space over
$\bF$ by letting $(f+g)(x) := f(x)+g(x)$ and $(k\cdot f)(x) := k\cdot
f(x)$ for all $x\in V$ and $k\in \bF$, and by setting the $K$-chain
$[x\mapsto 0]$ as the zero vector. It is worth noticing that
$\dim(K^V)= \dim(K)\cdot |V|$. If $K$ is equipped with a non-degenerate
$\sigma$-sesqui-bilinear form $\varphi$, we let $\langle,\rangle_\varphi:K^V\times
K^V\to \bF$ be such that for all $f,g\in K^V$,
\begin{align*}
\langle f,g \rangle_\varphi &:= \sum\limits_{x\in V} \varphi(f(x),g(x)).
\end{align*}

It is straightforward to verify that $\langle,\rangle_\varphi$ is a
non-degenerate $\sigma$-sesqui-bilinear form. (We will often write
$\langle,\rangle$ for convenience when the context is clear.)
Subspaces of $K^V$ are called \emph{$K$-chain groups on $V$}. A
$K$-chain group $L$ on $V$ is said \emph{lagrangian} if it is totally
isotropic and $\dim(L)=|V|$.

A \emph{simple isomorphism} from a $K$-chain group $L$ on $V$ to a
$K$-chain group $L'$ on $V'$ is a bijection $\mu:V\to V'$ such that
$L=\{f\circ \mu \mid f\in L'\}$ where $(f\circ \mu)(x) = f(\mu(x))$
for all $x\in V$. In this case we say that $L$ and $L'$ are
\emph{simply isomorphic}. 

From now on, we are only interested in $\bK_\sigma$-chain groups on
$V$.  Recall that $\bK_\sigma$ is the $2$-dimensional vector space
$\bF^2$ over $\bF$ equipped with the $\sigma$-sesqui-bilinear form
$\mathbf{b}_\sigma$. The following is a direct consequence of
definitions and Theorem \ref{thm:3.1}.

\begin{lem}\label{lem:3.1} If $L$ is a totally isotropic $\bK_\sigma$-chain
  group on $V$, then $\dim(L) \leq |V|$. If $L$ is lagrangian, then
  $L=L^\bot$.
\end{lem}

\begin{lem}\label{lem:3.2} Let $u,v\in \bK_\sigma$ and assume $u\ne 0$ is
  isotropic. If $\mathbf{b}_\sigma(u,v)=0$, then $v=c\cdot u$ for some
  $c\in \bF$.
\end{lem}

\begin{pf*}{Proof.} Since $\mathbf{b}_\sigma$ is non-degenerate, there
  exists $u'\in \bK_\sigma$ such that $\mathbf{b}_\sigma(u,u')\ne 0$. In this
  case, $\{u,u'\}$ is a basis for $\bK_\sigma$ (Property \ref{pt:3.1}). Hence,
  there exist $c,d\in \bF$ such that $v=c\cdot u + d\cdot
  u'$. Therefore,
  \begin{align*}
\mathbf{b}_\sigma (u,v) &= \frac{\sigma(c)}{\sigma(1)}\cdot
\mathbf{b}_\sigma (u,u) + \frac{\sigma(d)}{\sigma(1)}\cdot
\mathbf{b}_\sigma (u,u') = \frac{\sigma(d)}{\sigma(1)}\cdot
\mathbf{b}_\sigma (u,u').
  \end{align*}
Since $\mathbf{b}_\sigma (u,u') \ne 0$ and $\mathbf{b}_\sigma (u,v) =0$, we
have that $\sigma(d)=0$, \ie, $d=0$. \qed
\end{pf*}

We now introduce \emph{minors} for $\bK_\sigma$-chain groups on
$V$. If $f$ is a $\bK_\sigma$-chain on $V$, then $Sp(f):=\{x\in V\mid
f(x) \ne 0\}$.  If $L\subseteq \bK_\sigma^V$ and $X\subseteq V$, we
let $L_{\mid X} := \{\restriction{f}{ X}\mid f\in L\}$ and $L^{\mid X} :=
\{\restriction{f}{ X} \mid f\in L$ and $Sp(f)\subseteq X \}$. For $\alpha\in
\bK_\sigma^*$ and $X\subseteq V$, we let $L\minor{\alpha} X$ be the
$\bK_\sigma$-chain group
\begin{align*} L\minor{\alpha} X&:=\{\restriction{f}{ (V\setminus X)}\mid f\in L\ \textrm{and $\mathbf{b}_\sigma (
f(x), \alpha ) = 0$ for all}\ x\in X\}\end{align*} on $V\setminus X$.
  A pair $\{\alpha,\beta\}\subseteq \bK_\sigma^*$ is said
  \emph{minor-compatible} if $\mathbf{b}_\sigma( \alpha, \alpha) =
  \mathbf{b}_\sigma (\beta, \beta ) = 0$ and $\{\alpha,\beta\}$ forms a
  basis for $\bK_\sigma$. For a minor-compatible pair $\{\alpha,\beta\}$, a
  $\bK_\sigma$-chain group on $V\setminus (X\cup Y)$ of the form
  $L\minor{\alpha} X\minor{\beta}Y$ is called an
  \emph{$\alpha\beta$-minor} of $L$.

One easily verifies that $L\minor{\alpha} X\minor{\alpha} Y =
L\minor{\alpha} (X\cup Y)$, and $L\minor{\alpha} X\minor{\beta} Y=
L\minor{\beta}Y\minor{\alpha} X$. Hence, we have the following which
is already proved in \cite{OUM10} for a special case of
$\{\alpha,\beta\}$.

\begin{prop}\label{prop:3.1} Let $\{\alpha,\beta\}$ be
  minor-compatible. An $\alpha\beta$-minor of an $\alpha\beta$-minor
  of $L$ is an $\alpha\beta$-minor of $L$.
\end{prop}

We now prove that $\alpha\beta$-minors of lagrangian
$\bK_\sigma$-chain groups are also lagrangian. The proofs are the same
as in \cite{OUM10}. We include some of them that we expect can
convince the reader that the proofs are not different.

\begin{prop}\label{prop:3.2} Let $\{\alpha,\beta\}$ be
  minor-compatible. An $\alpha\beta$-minor of a totally
  isotropic $\bK_\sigma$-chain group $L$ on $V$ is totally isotropic.
\end{prop}

\begin{pf*}{Proof.} Let $L':=L\minor{\alpha} X\minor{\beta}Y$ be an $\alpha\beta$-minor
  of $L$ on $V':=V\setminus (X\cup Y)$. Let $f',g'\in L'$ and let
  $f,g\in L$ such that $f'=\restriction{f}{ V'}$ and $g'=\restriction{g}{V'}$. By Lemma
  \ref{lem:3.2}, for all $x\in X\cup Y$, $\mathbf{b}_\sigma(f(x),g(x)) =
  0$. Hence, $\sum\limits_{x\in V} \mathbf{b}_\sigma (f(x),g(x)) =
  \sum\limits_{x\in V'} \mathbf{b}_\sigma (f(x),g(x)) = \langle
  f',g'\rangle$. Therefore, $\langle f',g'\rangle = 0$. \qed
\end{pf*}

\begin{lem} \label{lem:3.3} Let $L$ be a $\bK_\sigma$-chain group on $V$ and
  $X\subseteq V$. Then, $\dim(L_{\mid X}) + \dim(L^{\mid (V\setminus X)}) = \dim(L)$
\end{lem}

\begin{pf*}{Proof.} Let $\varphi:L\to L_{\mid X}$ be the linear
  transformation that maps any $f\in L$ to $\restriction{f}{ X}$. We have
  clearly $L_{\mid X} = Im(\varphi)$. For any $f\in Ker(\varphi)$, we
  have $f(x)=0$ for all $x\in X$. Hence, $L^{\mid
    (V\setminus X)} = Ker(\varphi)$. This concludes the lemma. \qed
\end{pf*}

For any $x\in V$ and $\gamma \in \bK_\sigma^*$, we let $x^\gamma$ be
the $\bK_\sigma$-chain on $V$ such that
\begin{align*}
x^\gamma(z) &:=\begin{cases} \gamma & \textrm{if $z=x$},\\ 0 &
\textrm{otherwise}. \end{cases} 
\end{align*}

The following admits a similar proof as the one in \cite[Proposition
  3.6]{OUM10}. 

\begin{prop}\label{prop:3.3} Let $L$ be a $\bK_\sigma$-chain group on $V$,
  $x\in V$ and $\gamma \in \bK_\sigma^{*}$. Hence,
\begin{align*}
\dim(L\minor{\gamma} x) &= \begin{cases} \dim(L) & \textrm{if
    $x^\gamma\in L^\bot\setminus L$},\\ \dim(L) - 2 &
  \textrm{if $x^\gamma\in L\setminus L^\bot$}, \\ \dim(L) -
    1 & \textrm{otherwise}.
\end{cases} \end{align*}
\end{prop}

\begin{cor}\label{cor:3.1} Let $\{\alpha,\beta\}$ be
  minor-compatible. If $L$ is a totally isotropic $\bK_\sigma$-chain
  group on $V$ and $L'$ is an $\alpha\beta$-minor of $L$ on $V'$, then
  $|V'| - \dim(L') \leq |V| - \dim(L)$.
\end{cor}

\begin{pf*}{Proof.} By induction on $|V\setminus V'|$. Since $L$ is
  totally isotropic, for all $x\in V\setminus V'$, we cannot have
  neither $x^\alpha\in L\setminus L^\bot$ nor $x^\beta \in
  L\setminus L^\bot$. Hence, $\dim(L) - \dim(L\minor{\alpha} x)\in
  \{0,1\}$ and $\dim(L)-\dim(L\minor{\beta} x)\in \{0,1\}$ by
  Proposition \ref{prop:3.3}. Hence, if $|V\setminus V'|=1$, we are
  done.

  If $|V\setminus V'|>1$, let $x\in V\setminus V'$. Hence, $L'$ is
  an $\alpha\beta$-minor of $L\minor{\alpha} x$ or $L\minor{\beta}
  x$. By inductive hypothesis, $|V'| - \dim(L') \leq |V\setminus x| -
  \dim(L\minor{\alpha} x)$ or $|V'| - \dim(L') \leq |V\setminus x| -
  \dim(L\minor{\beta} x)$. And since, $|V\setminus x| -
  \dim(L\minor{\alpha} x) \leq |V| - \dim(L)$ and $|V\setminus x| -
  \dim(L\minor{\beta} x) \leq |V|-\dim(L)$, we are done. \qed
\end{pf*}

\begin{prop}\label{prop:3.4} Let $\{\alpha,\beta\}$ be
  minor-compatible. An $\alpha\beta$-minor of a lagrangian
  $\bK_\sigma$-chain group on $V$ is lagrangian.
\end{prop}

\begin{pf*}{Proof.} Let $L'$ be an $\alpha\beta$-minor of $L$ on
  $V'$. By Proposition \ref{prop:3.2}, $L'$ is totally isotropic, hence
  $\dim(L')\leq |V'|$. By Corollary \ref{cor:3.1}, $|V'|-\dim(L') \leq 0$
  since $\dim(L)=|V|$ ($L$ lagrangian). Hence, $\dim(L')\geq |V'|$. \qed
\end{pf*}

We now define the connectivity function for lagrangian
$\bK_\sigma$-chain groups. Let $L$ be a lagrangian $\bK_\sigma$-chain
group on $V$. For every $X\subseteq V$, we let
$\lambda_L(X):=|X|-\dim(L^{\mid X})$. Since $L^{\mid X}$ is totally
isotropic, $\dim(L^{\mid X}) \leq |X|$, and hence $\lambda_L(X)\geq
0$. 

\begin{prop}[\cite{OUM10}]\label{prop:3.5} Let $L$ be a lagrangian
  $\bK_\sigma$-chain group on $V$. Then, $\lambda_L$ is symmetric and
  submodular.
\end{prop}

The proof of Proposition \ref{prop:3.5} uses the fact that $2\cdot
\lambda_L(X) = \dim(L) - \dim(L^{\mid X}) - \dim(L^{\mid (V\setminus
  X)})$ and the following theorem by Tutte.

\begin{thm}[\cite{OUM10}]\label{thm:3.2} If $L$ is a $\bK_\sigma$-chain group on $V$ and
  $X\subseteq V$, then $(L_{\mid X})^\bot = (L^\bot)^{\mid X}$.
\end{thm}

The branch-width of a lagrangian $\bK_\sigma$-chain group $L$ on $V$, denoted
by $\bwd{L}$, is then defined as the $\lambda_L$-width of $V$. 

We can now state the well-quasi-ordering of lagrangian
$\bK_\sigma$-chain groups of bounded branch-width under
$\alpha\beta$-minor. Let us first enrich the $\alpha\beta$-minor to
labelled $\bK_\sigma$-chain groups on $V$. Let $(Q,\preceq)$ be a
well-quasi-ordered set. A \emph{$Q$-labelling} of a lagrangian
$\bK_\sigma$-chain group $L$ on $V$ is a mapping $\gamma_L:V\to Q$. A
\emph{$Q$-labelled} lagrangian $\bK_\sigma$-chain group on $V$ is a
couple $(L,\gamma_L)$ where $L$ is a lagrangian $\bK_\sigma$-chain
group on $V$ and $\gamma_L$ a $Q$-labelling of $L$. A $Q$-labelled
lagrangian $\bK_\sigma$-chain group $(L',\gamma_{L'})$ on $V'$ is an
\emph{$(\alpha\beta, Q)$-minor} of a $Q$-labelled lagrangian
$\bK_\sigma$-chain group $(L,\gamma_L)$ on $V$ if $L'$ is an
$\alpha\beta$-minor of $L$ and $\gamma_{L'}(x)\preceq \gamma_L(x)$ for
all $x\in V'$. $(L,\gamma_L)$ is \emph{simply isomorphic} to
$(L',\gamma_{L'})$ if there exists a simple isomorphism $\mu$ from $L$
to $L'$ and $\gamma_L = \gamma_{L'} \circ \mu $. The following is more
or less proved in \cite{OUM10}.

\begin{thm}\label{thm:3.3} Let $\bF$ be a finite field and $k$ a
  positive integer, and let $\{\alpha,\beta\}$ be
  minor-compatible. Let $(Q,\preceq)$ be a well-quasi-ordered set and
  let $(L_1,\gamma_{L_1}),(L_2,\gamma_{L_2}),\ldots$ be an infinite
  sequence of $Q$-labelled lagrangian $\bK_{\sigma_i}$-chain groups
  having branch-width at most $k$. Then, there exist $i<j$ such that
  $(L_i,\gamma_{L_i})$ is simply isomorphic to an
  $(\alpha\beta,Q)$-minor of $(L_j,\gamma_{L_j})$.
\end{thm}

Theorem \ref{thm:3.3} is proved in \cite{OUM10} for
$\alpha=\left(\begin{smallmatrix} 1\\ 0 \end{smallmatrix}
\right),\ \beta=\left(\begin{smallmatrix} 0\\ 1 \end{smallmatrix}
\right)$ and $\langle, \rangle_{\mathbf{b}_{\sigma_i}}$ being a (skew)
symmetric bilinear form. However, the proof uses only the axioms in
Properties \ref{pt:3.1} and \ref{pt:3.2}, and Theorem
\ref{thm:3.1}. The other necessary ingredients are Lemmas
\ref{lem:3.1}, \ref{lem:3.2} and \ref{lem:3.3}, Proposition
\ref{prop:3.3}, and Theorem \ref{thm:3.2}.  We refer to \cite{OUM10}
for the technical details. It is important that the reader keeps in
mind that even if $\mathbf{b}_\sigma$ is not a bilinear form, it
shares with the bilinear forms in \cite{OUM10} the necessary
properties for proving Theorem \ref{thm:3.3}.

\section{Representations of $\bK_\sigma$-Chain
  Groups by $(\sigma,\epsilon)$-Symmetric Matrices}\label{sec:4} 

In this section we will use Theorem \ref{thm:3.3} to obtain a similar
result for $(\sigma,\epsilon)$-symmetric matrices. We recall that we
use the Greek letter $\sigma$ for sesqui-morphisms, and if $\bF$ is a
field, then we let $\bK_\sigma$ be the $2$-dimensional vector space
$\bF^2$ over $\bF$ equipped with the $\sigma$-sesqui-bilinear form
$\mathbf{b}_\sigma$. We will associate with each
$(\sigma,\epsilon)$-symmetric matrix a lagrangian $\bK_\sigma$-chain
group. These matrices are called \emph{matrix representations}.  We
also need to relate $\alpha\beta$-minors of lagrangian
$\bK_\sigma$-chain groups to principal submatrices of their matrix
representations, and relate $\bF$-rank-width of
$(\sigma,\epsilon)$-symmetric matrices to branch-width of lagrangian
$\bK_\sigma$-chain groups. We follow similar steps as in \cite{OUM10}.

Let $\epsilon:V\to \{-1,+1\}$ be a function. We say that two
$\bK_\sigma$-chains $f$ and $g$ on $V$ are
\emph{$\epsilon$-supplementary} if, for all $x\in V$,
\begin{enumerate}
\item[(i)] $\mathbf{b}_\sigma(f(x),f(x)) = \mathbf{b}_\sigma(g(x),g(x)) =
  0$, 
\item[(ii)] $\mathbf{b}_\sigma(f(x),g(x)) = \epsilon(x)\cdot \sigma(1)$
  and
\item[(iii)] $\mathbf{b}_\sigma(g(x),f(x)) = -\epsilon(x)\cdot
  \sigma(1)^2$.
\end{enumerate}

For any $c\in \bF^*$, we let $c^*:=\left(\begin{smallmatrix}
  c\\ 0 \end{smallmatrix} \right)$, $c_*:=\left(\begin{smallmatrix}
  0\\ c \end{smallmatrix} \right)$, $\widetilde{c^*} :=
\left(\begin{smallmatrix} 0\\ \sigma(c^{-1}) \end{smallmatrix}
\right)$ and $\widetilde{c_*}:= \left(\begin{smallmatrix}
  -\sigma(1)\cdot \sigma(c)^{-1}\\ 0 \end{smallmatrix} \right)$. 

As a consequence of the following easy property, we get that for any
$\epsilon:V\to \{-1,+1\}$, we can construct $\epsilon$-supplementary
$\bK_\sigma$-chains on $V$. 

\begin{property}\label{pt:4.1} For any $c\in \bF^*$ and $\epsilon\in
  \{-1,+1\}$, we have  
  \begin{align*}
\begin{cases} 
\mathbf{b}_\sigma\left(\epsilon\cdot c^*, \widetilde{c^*}\right) & =
\epsilon\cdot \sigma(1)\\
\mathbf{b}_\sigma\left(\widetilde{c^*}, \epsilon\cdot c^* \right) &=
-\epsilon\cdot \sigma(1)^2  \end{cases} &\ \textrm{and}\ 
\begin{cases} 
\mathbf{b}_\sigma\left(\epsilon\cdot c_*, \widetilde{c_*}\right) & =
\epsilon\cdot \sigma(1)\\
\mathbf{b}_\sigma\left(\widetilde{c_*}, \epsilon\cdot c_* \right) &=
-\epsilon\cdot \sigma(1)^2  \end{cases} \end{align*}
\end{property}

The following associates with each $(\sigma,\epsilon)$-symmetric
$(V,V)$-matrix a lagrangian $\bK_\sigma$-chain group on $V$.  

\begin{prop}\label{prop:4.1} Let $M$ be a
  $(\sigma,\epsilon)$-symmetric  $(V,V)$-matrix over $\bF$, and let
  $f$ and $g$ be 
  $\epsilon$-supplementaty $\bK_\sigma$-chains on $V$. For every $x\in V$, we
  let $f_x$ be the $\bK_\sigma$-chain on $V$ such that, for all $y\in V$,
  \begin{align*}
    f_x(y) &:= \begin{cases}  m_{xx}\cdot f(x) + g(x) & \textrm{if
        $y=x$},\\ m_{xy} \cdot f(y) & \textrm{otherwise}. \end{cases}
  \end{align*}
  Then, the $\bK_\sigma$-chain group on $V$ denoted by $(M,f,g)$ and spanned
  by $\{f_x\mid x\in V\}$ is lagrangian.
\end{prop}

\begin{pf*}{Proof.} It is enough to prove that for all $x,y$, $\langle
  f_x,f_y \rangle = 0$ and the $f_x$'s  are linearly
  independent. 

  For all $x,y\in V$ and all $z\in V\setminus \{x,y\}$, $\mathbf{b}_\sigma(f_x(z),f_y(z)) =
  \mathbf{b}_\sigma(m_{xz} \cdot f(z),m_{yz}\cdot f(z)) = m_{xz}\cdot
    \sigma(m_{yz})\cdot \sigma(1)^{-1} \cdot \mathbf{b}_\sigma(f(z),f(z)) =
  0$. Hence for all $x,y\in V$, 
  \begin{align*}
\langle f_x,f_y \rangle &= \mathbf{b}_\sigma\left(f_x(x),f_y(x)\right) +
\mathbf{b}_\sigma\left(f_x(y),f_y(y)\right)\\ &=
\mathbf{b}_\sigma\left(m_{xx}\cdot f(x) +
g(x), m_{yx}\cdot f(x)\right) + \mathbf{b}_\sigma\left(m_{xy}\cdot f(y),
m_{yy}\cdot f(y) + g(y)\right) \\ & = \sigma(m_{yx})\cdot \sigma(1)^{-1}\cdot \mathbf{b}_\sigma\left(g(x),f(x)\right) +
m_{xy}\cdot \mathbf{b}_\sigma\left(f(y),g(y)\right) \\ &= \sigma(1)\cdot
\left(\epsilon(y)\cdot m_{xy} -\epsilon(x)\cdot
\sigma(m_{yx}) \right) \\ & = 0.
  \end{align*}

  It remains to prove that the $f_x$'s are linearly
  independent. Assume there exist constants $c_x$ such that
  $\sum\limits_{x\in V} c_x\cdot f_x = 0$. Hence, for all $y\in V$,
  $\mathbf{b}_\sigma\left(f(y), \sum\limits_{x\in V} c_x\cdot
  f_x(y)\right) = 0$. But for all $x\in V$ and all $y\in V\setminus x$, $\mathbf{b}_\sigma\left(f(y), c_x\cdot
  f_x(y)\right) = 0$. Hence, for all $y\in V$, $\mathbf{b}_\sigma\left(f(y), \sum\limits_{x\in V} c_x\cdot
  f_x(y)\right) = \mathbf{b}_\sigma(f(y),c_y\cdot f_y(y)) =
  \epsilon(y)\cdot \sigma(c_y)$,
  \ie, $\sigma(c_y)= 0$. Hence, we conclude that $c_y=0$ for all $y\in
  V$, \ie, the $f_x$'s are linearly independent. \qed
\end{pf*}

If a lagrangian $\bK_\sigma$-chain group $L$ is simply isomorphic to
$(M,f,g)$, we call $(M,f,g)$ a \emph{matrix representation of
  $L$}. One easily verifies from the definition of $(M,f,g)$, that for
all non zero $\bK_\sigma$-chains $h\in (M,f,g)$, we do not have
$\mathbf{b}_\sigma(h(x),f(x)) = 0$ for all $x\in V$. We now make
precise this property.

A $\bK_\sigma$-chain $f$ on $V$ is called an \emph{eulerian chain} of a
lagrangian $\bK_\sigma$-chain group $L$ on $V$ if:
\begin{enumerate}
\item[(i)] for all $x\in V$, $f(x)\ne 0$ and
  $\mathbf{b}_\sigma(f(x),f(x))=0$, and   
\item[(ii)] there is no non-zero $\bK_\sigma$-chain $h$ in $L$ such that
  $\mathbf{b}_\sigma(h(x),f(x))=0$ for all $x\in V$. 
\end{enumerate}

The proof of the following is the same as in \cite{OUM10}. 

\begin{prop}[\cite{OUM10}]\label{prop:4.2} Every lagrangian $\bK_\sigma$-chain group on $V$
  has an eulerian chain.
\end{prop}

\begin{pf*}{Proof.} By induction on the size of $V$.  We let
  $\alpha:=c^*$ and $\beta:= \widetilde{c^*}$ for some $c\in
  \bF^*$. Let $L$ be a lagrangian $\bK_\sigma$-chain group on $V$. If
  $V=\{x\}$, then $\dim(L)=1$, hence either $x^\alpha$ or $x^\beta$ is
  an eulerian chain.

  Assume $|V|>1$ and let $V':=V\setminus x$ for some $x\in V$. Hence,
  both $L\minor{\alpha} x$ and $L\minor{\beta} x$ are lagrangian. By
  inductive hypothesis, there exist $f'$ and $g'$ such that $f'$
  (resp. $g'$) is an eulerian chain of
  $L\minor{\alpha} x$ (resp. $L\minor{\beta} x$).

  Let $f$ and $g$ be $\bK_\sigma$-chains on $V$ such that
  $f(x)=\alpha$, $g(x)=\beta$, and $f'=\restriction{f}{ V'}$ and
  $g'=\restriction{g}{ V'}$. We claim that either $f$ or $g$ is an
  eulerian chain of $L$. Otherwise, there exist non-zero
  $\bK_\sigma$-chains $h$ and $h'$ in $L$ such that
  $\mathbf{b}_\sigma(h(x),f(x))=0$ and
  $\mathbf{b}_\sigma(h'(x),g(x))=0$ for all $x\in V$. Hence, we have
  $\mathbf{b}_\sigma(\restriction{h}{ V'}(x),f'(x))=0$ and
  $\mathbf{b}_\sigma(\restriction{h'}{ V'}(x),g'(x))=0$ for all $x\in
  V'$. Therefore, $\restriction{h}{ V'} = \restriction{h'}{ V'} = 0$,
  otherwise there is a contradiction because $\restriction{h}{ V'}\in
  L\minor{\alpha} x$ and $\restriction{h'}{ V'}\in L\minor{\beta} x$
  by construction of $f$ and $g$. Thus, $h(x)\ne 0$ and $h'(x)\ne 0$,
  and $\langle h,h'\rangle = \mathbf{b}_\sigma(h(x),h'(x))$. By Lemma
  \ref{lem:3.2}, we have $h(x)=d\cdot \alpha$ and $h'(x)=d'\cdot
  \beta$ for some $d,d'\in \bF^*$. Hence, $\langle h,h'\rangle =
  d\cdot \sigma(d')\ne 0$, which contradicts the totally isotropy of
  $L$.\qed
\end{pf*}

The next proposition shows how to construct a matrix representation of
a lagrangian $\bK_\sigma$-chain group. 

\begin{prop}\label{prop:4.3} Let $L$ be a lagrangian $\bK_\sigma$-chain group
  on $V$. Let $\epsilon:V\to \{-1,+1\}$, and let $f$ and $g$ be
  $\epsilon$-supplementary with $f$ being an eulerian chain
  of $L$. For every $x\in V$, there exists a unique $\bK_\sigma$-chain
  $f_x\in L$ such that
  \begin{enumerate}
  \item[(i)] $\mathbf{b}_\sigma(f(y),f_x(y)) = 0$ for all $y\in
    V\setminus x$, 
  \item[(ii)] $\mathbf{b}_\sigma(f(x),f_x(x)) = \epsilon(x)\cdot \sigma(1)$.
  \end{enumerate}
  Moreover, $\{f_x\mid x\in V\}$ is a basis for $L$. If we let $M$ be
  the $(V,V)$-matrix such that $m_{xy}:=\mathbf{b}_\sigma(f_x(y),g(y))\cdot
  \sigma(1)^{-1}\cdot \epsilon(y)$, then $M$ is
  $(\sigma,\epsilon)$-symmetric and $(M,f,g)$ is a matrix
  representation of $L$.
\end{prop}

\begin{pf*}{Proof.} The proof is the same as the one in
  \cite{OUM10}. We first prove that $\bK_\sigma$-chains verifying 
  statements (i) and (ii) exist. For every $x\in V$, let $g_x$ be the
  $\bK_\sigma$-chain on $V$ such that $g_x(x)=f(x)$ and $g_x(y)=0$ for all
  $y\in V\setminus x$. We let $W$ be the $\bK_\sigma$-chain group spanned by
  $\{g_x\mid x\in V\}$. The dimension of $W$ is clearly $|V|$. Let
  $L+W=\{h+h'\mid h\in L,\ h'\in W\}$. We have $L\cap W=\{0\}$ because
  $f$ is eulerian to $L$. Hence, $\dim(L+W)= 2\cdot |V|$,
  \ie, $\bK_\sigma^V=L+W$.  For each $x\in V$, let $h_x\in \bK_\sigma^V$ such that
  $h_x(x)=g(x)$ and $h_x(y)=0$ for all $y\in V\setminus x$. Hence,
  there exist $f_x\in L$ and $g'_x\in W$ such that $h_x=f_x+g'_x$. We
  now prove that these $f_x$'s verify statements (i) and (ii). Let
  $g'_x=\sum\limits_{z\in V} c_z\cdot g_z$. For all $x\in V$ and all
  $y\in V\setminus x$,
  \begin{align*}
    \mathbf{b}_\sigma(f(x),f_x(x)) &= \mathbf{b}_\sigma(f(x), h_x(x)-g'_x(x))
    \\ & = \mathbf{b}_\sigma(f(x), h_x(x)) -  \mathbf{b}_\sigma(f(x), 
    g'_x(x)) \\ &= \mathbf{b}_\sigma(f(x),g(x)) -
    \mathbf{b}_\sigma(f(x),c_x\cdot f(x)) \\ & = \epsilon(x)\cdot
    \sigma(1) \\ \intertext{and}
    \mathbf{b}_\sigma(f(y),f_x(y)) &= \mathbf{b}_\sigma(f(y),h_x(y)) -
    \mathbf{b}_\sigma(f(y),c_y\cdot g_y(y)) \\ &= \mathbf{b}_\sigma(f(y),0) -
    \mathbf{b}_\sigma(f(y),c_y\cdot f(y)) = 0.
  \end{align*}
  We now prove that each $f_x$ is unique. Assume there exist $f_x$'s
  and $f'_x$'s verifying statements (i) and (ii). For each $x\in V$, we have
  $\mathbf{b}_\sigma(f(x),f_x(x)-g(x)) = \mathbf{b}_\sigma(f(x),f_x(x)) -
  \mathbf{b}_\sigma(f(x),g(x)) = 0$. Similarly,
  $\mathbf{b}_\sigma(f(x),f'_x(x)-g(x))=0$. Hence, by Lemma
  \ref{lem:3.2}, $f_x(x)= c\cdot f(x)+g(x)$ and $f'_x(x)=c'\cdot
  f(x)+g(x)$ for $c,c'\in \bF^*$. We let $h'_x=f_x-f'_x$ which belongs
  to $L$. Therefore, for all
  $z\in V$, we have $\mathbf{b}_\sigma(f(z),h'_x(z))= 0$. And since $f$
  is eulerian to $L$, we have $h'_x=0$, \ie, $f_x=f'_x$. 

  By using the same technique as in the proof of Proposition
  \ref{prop:4.1}, one easily proves that $\{f_x\mid x\in V\}$ is
  linearly independent.  It remains to prove that $M:=(m_{xy})_{x,y\in
    V}$ with
  $m_{xy}=\mathbf{b}_\sigma(f_x(y),g(y))\cdot \sigma(1)^{-1}\cdot
  \epsilon(y)$ is $(\sigma,\epsilon)$-symmetric and
  $L=(M,f,g)$.

  We recall that $f(x)$ is isotropic for all $x\in V$. By statement
  (i) and Lemma \ref{lem:3.2}, for all $x\in V$ and all $y\in
  V\setminus x$, we have $f_x(y)=c_{xy}\cdot f(y)$ for some
  $c_{xy}\in \bF$. Hence, $m_{xy} = c_{xy}$. Similarly, we have
  $f_x(x) = c_{xx}\cdot f(x) + g(x)$ for some $c_{xx}\in \bF$, \ie,
  $m_{xx} = c_{xx}$.  It is thus clear that $L=(M,f,g)$. We now show
  that $M$ is $(\sigma,\epsilon)$-symmetric. Since $L$ is isotropic,
  we have for all $x,y\in V$, $\langle f_x,f_y\rangle =
  \mathbf{b}_\sigma(f_x(x),f_y(x)) +
  \mathbf{b}_\sigma(f_x(y),f_y(y))=0$. But,
  \begin{align*}
\mathbf{b}_\sigma(f_x(x),f_y(x)) +
  \mathbf{b}_\sigma(f_x(y),f_y(y)) &= \mathbf{b}_\sigma(m_{xx}\cdot
  f(x)+g(x), m_{yx}\cdot f(x)) +\\ & \qquad \qquad\qquad \mathbf{b}_\sigma(m_{xy}\cdot f(y),
  m_{yy}\cdot f(y) + g(y)) \\ &= \sigma(m_{yx})\cdot
  \sigma(1)^{-1}\cdot \mathbf{b}_\sigma(g(x),f(x)) + m_{xy}\cdot
  \mathbf{b}_\sigma(f(y),g(y)) \\ &= \sigma(1)\cdot (\epsilon(y)\cdot
  m_{xy} - \epsilon(x)\cdot \sigma(m_{yx}))
  \end{align*}
  Hence, $\epsilon(y)\cdot m_{xy} = \epsilon(x)\cdot
  \sigma(m_{yx})$. \qed
\end{pf*}

From Proposition \ref{prop:4.1} (resp. \ref{prop:4.3}), to every
every $(\sigma,\epsilon)$-symmetric $(V,V)$-matrix (resp. lagrangian
$\bK_\sigma$-chain group on $V$) one can associate a lagrangian
$\bK_\sigma$-chain group on $V$ (resp. a $(\sigma,\epsilon)$-symmetric
$(V,V)$-matrix). The next theorem relates the branch-width of a
lagrangian $\bK_\sigma$-chain group on $V$ to the $\bF$-rank-width of
its matrix-representations. Its proof is present in \cite{OUM10}, but
we give it for completeness.

\begin{thm}[\cite{OUM10}]\label{thm:4.1} Let $(M,f,g)$ be a matrix representation of
  a lagrangian $\bK_\sigma$-chain group $L$ on $V$. For every $X\subseteq V$,
  we have $\fcutrk_M(X)=\lambda_L(X)$. 
\end{thm}

\begin{pf*}{Proof.} We let $\{f_x\mid x\in V\}$ be the basis of $L$
  given in Proposition \ref{prop:4.1}. Let
  $A:=\matind{M}{X}{V\setminus X}$. It is well-known in linear
  algebra that $\rk(A)=\rk(A^t)=|X|-n(A^t)$ where $n(A^t)$ is
  $\dim\left(\{p\in \bF^X\mid A^t\cdot p=0\}\right) = \dim\left(\{p\in
  \bF^X\mid p^t\cdot A =0\}\right)$. Let $\varphi:\bF^V\to L$ be such
  that $\varphi(p):= \sum\limits_{x\in V} p(x)\cdot f_x$. It is clear
  that $\varphi$ is a linear transformation and is therefore an
  isomorphism. Hence, 
  \begin{align*}
    \dim(L^{\mid X}) &= \dim\left(\{h\in L\mid Sp(h)\subseteq
    X\}\right)\\ &= \dim\left(\varphi^{-1}\left(\{h\in L\mid
    Sp(h)\subseteq X\}\right) \right)\\ &= \dim\left( \{p\in \bF^V\mid
    \sum\limits_{x\in V} p(x)\cdot f_x(y) = 0\ \textrm{for all}\ y\in
    V\setminus X\} \right).\end{align*} Now, let $p\in \bF^V$ such
  that $\restriction{\varphi(p)}{ X} \in L^{\mid X}$. Then, for all
  $y\in V\setminus X$, $\varphi(p)(y) =0$, \ie,
  $\mathbf{b}_\sigma(f(y),\varphi(p)(y))= 0$. But, $\varphi(p)(y)=
  \sum\limits_{x\in V} p(x)\cdot f_x(y)$. And, since
  $\mathbf{b}_\sigma(f(y), f_x(y))=0$ for all $x\ne y$, we have
  $\mathbf{b}_\sigma(f(y), \varphi(p)(y)) = \mathbf{b}_\sigma(f(y),
  p(y)\cdot f_y(y)) = \sigma(p(y))\cdot \epsilon(y)$, \ie,
  $p(y)=0$. Hence,
  \begin{align*}
    \dim(L^{\mid X}) &= \dim\left( \{p\in \bF^X\mid \sum\limits_{x\in
      X} p(x)\cdot m_{xy} = 0\ \textrm{for all}\ y\in V\setminus X\}
    \right) \\ & = \dim\left (\{p\in \bF^X\mid p^t\cdot A= 0\}\right)
    \\ &= n(A^t) \end{align*} Since, $\lambda_L(X)=|X|-\dim(L^{\mid X})$,
  we can conclude that $\fcutrk_M(X)=\lambda_L(X)$.  \qed
\end{pf*}

It remains now to relate $\alpha\beta$-minors of lagrangian
$\bK_\sigma$-chain groups to principal submatrices of their matrix
representations. For doing so, we need to prove some technical lemmas.
For $X \subseteq V$, we let $P_X$ and $I_X$ be the non-singular
diagonal $(V,V)$-matrices where
\begin{align*}
 P_X[x,x] &:= \begin{cases} \sigma(-1) & \textrm{if $x\in X$},\\ 1 &
  \textrm{otherwise}, \end{cases} & \quad \textrm{and} & \quad I_X[x,x]
 := \begin{cases} -1 & \textrm{if $x\in X$},\\ 1 &
  \textrm{otherwise.} \end{cases}
\end{align*}

If $M$ is a matrix of the form $\left(\begin{smallmatrix} \alpha & \beta \\ \gamma &
  \delta \end{smallmatrix}\right)$ where $\alpha:=M[X]$ is non-singular, the
\emph{principal pivot transform} of $M$ at $X$, denoted by $M*X$, is
the matrix
\begin{align*} \centering 
  \begin{pmatrix} \alpha^{-1} & \alpha^{-1}\cdot \beta \\ -\gamma \cdot
    \alpha^{-1} & \ \ M/\alpha
  \end{pmatrix}.
\end{align*}

The principal pivot transform was introduced by Tucker \cite{TUC60} in
an attempt to understand the linear algebraic structure of the
\emph{simplex method} by Dantzig. It appeared to have wide
applicability in many domains; without being exhaustive we can cite
linear algebra \cite{TSAT00}, graph theory \cite{BOU88} and biology
\cite{BH10}.

\begin{prop}\label{prop:4.4} Let $(M,f,g)$ be a matrix representation
  of a lagrangian $\bK_\sigma$-chain group $L$ on $V$. Let $X\subseteq V$
  such that $M[X]$ is non-singular. Let $f'$ and $g'$ be $\bK_\sigma$-chains
  on $V$ such that, for all $x\in V$,
  \begin{align*}
    f'(x) &:= \begin{cases} f(x) & \textrm{if $x\notin
        X$},\\ g(x) &
      \textrm{otherwise}, \end{cases} & \quad \textrm{and} & \quad g'(x) &:= \begin{cases}
     g(x) &
      \textrm{if $x\notin X$},\\ \sigma(-1)\cdot f(x) &
      \textrm{otherwise}. \end{cases} \end{align*}
  Then, $(P_X\cdot (M*X),f',g')$ is a matrix representation of $L$.  
\end{prop}

\begin{pf*}{Proof.} Let $\epsilon$ be such that $M$ is
  $(\sigma,\epsilon)$-symmetric, \ie, $f$ and $g$ are
  $\epsilon$-supplementary. Let us first show that $f'$ and $g'$ are
  $\epsilon$-supplementary. Since for all $x\notin X$, we have
  $f'(x)=f(x)$ and $g'(x)=g(x)$, we need to verify the properties of
  $\epsilon$-supplementary for the $x\in X$. For each
  $x\in X$, we have:
  \begin{align*}
\mathbf{b}_\sigma(f'(x),f'(x)) &= 
\mathbf{b}_\sigma(g(x),g(x)) =
0\\\\ \mathbf{b}_\sigma(g'(x),g'(x)) &=
\mathbf{b}_\sigma(\sigma(-1)\cdot f(x),\sigma(-1)\cdot f(x))\\ &=\mathbf{b}_\sigma(f(x),f(x)) =
0\\\\ \mathbf{b}_\sigma(f'(x),g'(x)) &=
\mathbf{b}_\sigma(g(x),\sigma(-1)\cdot f(x)) \\ &=\frac{-1}{\sigma(1)}
\cdot \mathbf{b}_\sigma(g(x),f(x)) = \epsilon(x)\cdot
\sigma(1) \\ \\ \mathbf{b}_\sigma(g'(x),f'(x)) &=
\mathbf{b}_\sigma(\sigma(-1)\cdot f(x),g(x)) \\ &=
-\sigma(1) \cdot
\mathbf{b}_\sigma(f(x),g(x)) = -\epsilon(x)\cdot \sigma(1)^2
  \end{align*}
 Hence, $f'$ and $g'$ are $\epsilon$-supplementary. It remains to show
 that $f'$ is eulerian to $L$.  For each $x\in V$, we let $f_x$ be the $\bK_\sigma$-chain on $V$ such
 that
  \begin{align*}
    f_x(y) &:= \begin{cases} m_{xy}\cdot f(y) & \textrm{if $y\ne
        x$},\\ m_{xx}\cdot f(x) + g(x) &
      \textrm{otherwise} \end{cases}
  \end{align*}
  By Propositions \ref{prop:4.1} and \ref{prop:4.3} the set $\{f_x\mid
  x\in V\}$ is a basis for $L$. Let $h\in L$ such that
  $\mathbf{b}_\sigma(h(y),f'(y))=0$ for all $y\in V$. Let
  $h=\sum\limits_{z\in V} c_z\cdot f_z$. For all $y\notin X$, we have
  \begin{align*}
\mathbf{b}_\sigma(h(y),f'(y)) &= 
\mathbf{b}_\sigma\left(\sum\limits_{z\in V}\left(c_z\cdot m_{zy}\cdot
f(y)\right)+c_y\cdot g(y), f(y)\right) \\ &= \mathbf{b}_\sigma(c_y\cdot
g(y), f(y)) \\ &= - c_y\cdot 
\epsilon(y)\cdot \sigma(1)^2.
  \end{align*}
Hence, $c_y=0$ for all $y\notin X$. If $y\in X$, then
  \begin{align*}
\mathbf{b}_\sigma(h(y),f'(y)) & =
\mathbf{b}_\sigma\left(\sum\limits_{z\in X}\left(c_z\cdot m_{zy}\cdot
f(y)\right) + c_y\cdot g(y), g(y)\right) \\ &=
\sum\limits_{z\in X} \left(c_z \cdot m_{z y}\cdot \mathbf{b}_\sigma(f(y),g(y))\right)
\\ & =\sigma(1)\cdot \epsilon(y)\cdot \sum\limits_{z\in X} c_z
\cdot m_{z y}. \end{align*} And for $\mathbf{b}_\sigma(h(y),f'(y))$ to
  being $0$, we must have $\sum\limits_{z\in X} \left(c_z \cdot m_{z
    y}\right) = 0$.  But, since $M[X]$ is non-singular, we have
  $\sum\limits_{z\in X} \left(c_z \cdot m_{z y}\right) = 0$ for all
  $y\in X$ if and only if $c_z=0$ for all $z\in X$. Therefore, we have
  $h=0$, \ie, $f'$ is eulerian.

By Proposition \ref{prop:4.3} there exists a unique matrix $M'$ such
that $L=(M',f',g')$. We will show that $M'=P_X\cdot (M*X)$. Assume
$M=\left( \begin{smallmatrix} \alpha & \beta \\ \gamma &
  \delta \end{smallmatrix} \right)$ with $\alpha := M[X]$. Let $I_f$
and $I_{\bar{f}}$ be respectively $(X,X)$ and $(V\setminus
X,V\setminus X)$-diagonal matrices with diagonal entries being the
$f(x)$'s. We define similarly, $I_g$ and $I_{\bar{g}}$, but diagonal
entries are $g(x)$'s. We let $A$ be the $(V,V)$-matrix, where $a_{xy}
:= f_x(y)$. Hence,
 \begin{align*}
    A &= \begin{pmatrix} \alpha \cdot I_f + I_g & \beta \cdot
    I_{\bar{f}} \\ \gamma \cdot I_f & \delta \cdot I_{\bar{f}} +
    I_{\bar{g}} \end{pmatrix}.
 \end{align*}
 The row space of $A$ is exactly $L$.  Let $B$ be the non-singular
 $(V,V)$-matrix 
  \begin{align*}
  \begin{pmatrix}\alpha^{-1} & 0 \\ 
    -\gamma \cdot \alpha^{-1} &\quad 
     I \end{pmatrix}. 
  \end{align*}  Therefore, 
  \begin{align*}
     B\cdot A & = \begin{pmatrix}\alpha^{-1}\cdot
      I_g + I_f &
      \alpha^{-1} \cdot \beta \cdot I_{\bar{f}} \\ -\gamma \cdot \alpha^{-1} \cdot I_g &(\delta-\gamma\cdot
      \alpha^{-1} \cdot \beta) \cdot I_{\bar{f}} +
       I_{\bar{g}} \end{pmatrix}.
  \end{align*}
    
 Let $A':=P_X\cdot B\cdot A$, and for each $x\in V$, let $f'_x$ be the
 $\bK_\sigma$-chain on $V$ with $f'_x(y):=a'_{xy}$.  From above, we
 have that $\{f'_x\mid x\in V\}$ is a basis for $L$. Let $C:=P_X\cdot
 (M*X)$. Then, for every $x,y\in V$, we have
 \begin{align*}
     f'_x(y) &= \begin{cases}
      c_{xy}\cdot f(y) &
      \textrm{if $y\ne x$ and $y\notin X$},\\ c_{xy}\cdot g(y) &
      \textrm{if $y\ne x$ and $y\in X$}, \\ c_{xx}\cdot f(x) + g(x) & \textrm{if $y=x\notin X$},
      \\ c_{xx}\cdot g(x) + \sigma(-1)\cdot f(x) & \textrm{if $y=x\in
        X$}. \end{cases} \\ 
\intertext{Hence,}
   \mathbf{b}_\sigma(f'(y),f'_x(y)) &= \begin{cases} \mathbf{b}_\sigma(f(y),c_{xy}\cdot f(y)) & \textrm{if $y\ne x$
       and $y\notin X$}, \\ \mathbf{b}_\sigma(g(y),c_{xy}\cdot g(y)) &
     \textrm{if $y\ne x$ and $y\in X$},\\ \mathbf{b}_\sigma(f(x),c_{xx}\cdot f(x)
     + g(x)) & \textrm{if $y=x\notin X$},\\
     \mathbf{b}_\sigma(g(x),c_{xx}\cdot g(x) + \sigma(-1)\cdot f(x)) &
     \textrm{if $y=x\in X$}.
     \end{cases}
 \end{align*} Hence, for all $x\in V$ and all $y\in V\setminus x$, we have
   $\mathbf{b}_\sigma(f'(x),f'_x(x)) =\epsilon(x)\cdot \sigma(1)$ and
 $\mathbf{b}_\sigma(f'(y),f'_x(y)) =0$.  Therefore, by Propositions
 \ref{prop:4.1} and \ref{prop:4.3} $\{f'_x\mid x\in V\}$ is the basis
 associated with $(M',f',g')$ and $M'=C=P_X\cdot (M*X)$. \qed

\end{pf*}


\begin{prop}\label{prop:4.5} Let $(M,f,g)$ be a matrix representation of
  a lagrangian $\bK_\sigma$-chain group $L$ on $V$ and let $Z\subseteq V$. Let
  $f'$ and $g'$ be $\bK_\sigma$-chains on $V$ such that 
  \begin{align*}
    f'(x)&:=\begin{cases} -f(x)& \textrm{if $x\in Z$},\\ f(x) &
    \textrm{otherwise}, \end{cases}  & \quad \textrm{and} & \quad g'(x) &:=\begin{cases} -g(x) &
    \textrm{if $x\in Z$}, \\ g(x) & \textrm{otherwise}. \end{cases}
  \end{align*}
  Then, $(I_Z\cdot M,f,g')$ and $(M\cdot I_Z,f',g)$ are matrix
  representations of $L$.
\end{prop}

\begin{pf*}{Proof.} Let $\epsilon:V\to \{+1,-1\}$ be such that $M$ is
  $(\sigma,\epsilon)$-symmetric, \ie, $f$ and $g$ are
  $\epsilon$-supplementary. Let $\{f_x\mid x\in V\}$ be the basis of
  $L$ associated with $f$ and $g$ by Proposition \ref{prop:4.1}.  One easily verifies that $f'$ and $g$, and $f$ and
  $g'$ are $\epsilon'$-supplementary with $\epsilon'(x)=-\epsilon(x)$
  if $x\in Z$, otherwise $\epsilon'(x)=\epsilon(x)$.  Moreover, $f'$
  is eulerian (because $f$ is eulerian). By Proposition
  \ref{prop:4.3}, there exist unique $f'_x$'s and $f''_x$'s such that
  $(M',f',g)$ and $(M'',f,g')$ are matrix representations of $L$ with
  $m'_{xy} := \mathbf{b}_\sigma(f'_x(y),g'(y))\cdot
  \sigma(1)^{-1}\cdot \epsilon'(y)$ and
  $m''_{xy}:=\mathbf{b}_\sigma(f''_x(y),g(y))\cdot \sigma(1)^{-1}\cdot
  \epsilon'(y)$.

  One easily checks that $\{-f_x\mid x\in Z\} \cup \{f_x\mid x\in
  V\setminus Z\}$ is the basis of $L$ associated with $f$ and $g'$ by
  Proposition \ref{prop:4.3}. It remains to prove that $M'=M\cdot
  I_Z$. If $x,y\in Z$, then
  $m'_{xy}=\mathbf{b}_\sigma(-f_x(y),-g(y))\cdot
  (-\epsilon(y))\cdot\sigma(1)^{-1}=-m_{xy}$. If $x\in Z$ and $y\notin
  Z$, then $m'_{xy}=\mathbf{b}_\sigma(-f_x(y),g(y))\cdot
  \epsilon(y)\cdot\sigma(1)^{-1}=-m_{xy}$. If $x,y\notin Z$, then
  $m'_{xy}=\mathbf{b}_\sigma(f_x(y),g(y))\cdot
  \epsilon(y)\cdot\sigma(1)^{-1}=m_{xy}$. And finally if $x\notin Z$
  and $y\in Z$, $m'_{xy}=\mathbf{b}_\sigma(f_x(y),-g(y))\cdot
  (-\epsilon(y))\cdot\sigma(1)^{-1}=m_{xy}$. Therefore, $M'=I_Z\cdot
  M$.

  It is straightforward to check that $\{f_x\mid x\in V\}$ is the
  basis of $L$ associated with $f'$ and $g$ by Proposition
  \ref{prop:4.3}. Then, $f''_x=f_x$. Let $x\in V$. We have clearly
  that $m''_{xy}=m_{xy}$ for all $y\in V\setminus Z$. Let now $y\in
  Z$. Hence, $m''_{xy}=-\mathbf{b}_\sigma(f_x(y),g(y))\cdot
  \epsilon(y)\cdot \sigma(1)^{-1}=-m_{xy}$. Hence, $M''=M\cdot
  I_Z$. \qed
\end{pf*}

A pair $(p,q)$ of non-zero scalars in $\bF$ is said
$\sigma$-compatible if $p^{-1}=\sigma(q)\cdot \sigma(1)^{-1}$
(equivalently $q^{-1}=\sigma(p)\cdot \sigma(1)^{-1}$). That means that
$(q,p)$ is also $\sigma$-compatible. It is worth noticing that if
$(p,q)$ is $\sigma$-compatible, then $(p^{-1},q^{-1})$ is also
$\sigma$-compatible. A pair $(P,Q)$ of non-singular diagonal
$(V,V)$-matrices is said $\sigma$-compatible if $(p_{xx},q_{xx})$ is
$\sigma$-compatible for all $x\in V$. For instance the pair
$(P_X,P_X^{-1})$ is $\sigma$-compatible.

\begin{prop}\label{prop:4.7}  Let $(M,f,g)$ be a matrix representation of
  a lagrangian $\bK_\sigma$-chain group $L$ on $V$ and let $(P,Q)$ be a
  $\sigma$-compatible pair of diagonal $(V,V)$-matrices. Let $f'$ and $g'$ be
  $\bK_\sigma$-chains on $V$ such that for all $x\in V$, $f'(x):=q_{xx}\cdot
  f(x)$ and $g'(x):=p_{xx}\cdot g(x)$. Then, $(P\cdot M\cdot
  Q^{-1},f',g')$ is a matrix representation of $L$.
\end{prop}

\begin{pf*}{Proof.} Let $\epsilon:V\to \{+1,-1\}$ such that $M$ is 
  $(\sigma,\epsilon)$-symmetric, \ie, $f$ and $g$ are
  $\epsilon$-supplementary. It is a straightforward computation to
  check that $f'$ and $g'$ are $\epsilon$-supplementary
  $\bK_\sigma$-chains on $V$. Moreover, $f'$ is eulerian to $L$
  (because $f$ is). By Proposition \ref{prop:4.3}, there exists a
  unique basis $\{f'_x\mid x\in V\}$ of $L$ such that $(M',f',g')$ is
  a matrix representation of $L$ with $m'_{xy}:=
  \mathbf{b}_\sigma(f'_x(y),g'(y)) \cdot \epsilon(y)\cdot
  \sigma(1)^{-1}$. Let $\{f_x\mid x\in V\}$ be the basis of $L$
  associated with $f$ and $g$ by Proposition \ref{prop:4.1}.

  For each $x\in V$, we clearly have
  $\mathbf{b}_\sigma(f'(y),p_{xx}\cdot f_x(y))=q_{yy}\cdot
  q_{xx}^{-1}\cdot \mathbf{b}_\sigma(f(y),f_x(y))$ for all
  $x,y\in V$. Therefore, for all $x\in V$ and all $y\in V\setminus x$,
  we have 
  \begin{align*}
    \mathbf{b}_\sigma(f'(x),p_{xx}\cdot f_x(x)) &= \epsilon(x)\cdot
    \sigma(1),\\
    \mathbf{b}_\sigma(f'(y),p_{xx}\cdot f_x(y)) &= 0.
\end{align*}
  Hence, by Proposition \ref{prop:4.3} $f'_x=p_{xx}\cdot f_x$. Then, for each $x,y\in
  V$, we have
  \begin{align*}
    m'_{xy}&=\mathbf{b}_\sigma(p_{xx}\cdot f_x(y), p_{yy}\cdot
    g(y))\cdot \epsilon(y)\cdot \sigma(1)^{-1}\\ &= p_{xx}\cdot
    \sigma(p_{yy})\cdot \sigma(1)^{-1} \cdot \left(\mathbf{b}_\sigma(f_x(y),
    g(y))\cdot \epsilon(y)\cdot \sigma(1)^{-1}\right) = p_{xx}\cdot
    q_{yy}^{-1}\cdot m_{xy}.
  \end{align*}
  Hence, $(P\cdot M\cdot Q^{-1},f',g')$ is a matrix representation of
  $L$. \qed
\end{pf*}

We call $(M,f,g)$ a \emph{special matrix representation} of a
lagrangian $\bK_\sigma$-chain group $L$ on $V$ if $f(x),g(x)\in
\{c^*,c_*\mid c\in \bF^*\}$ for all $x\in V$. A special case of the
following is proved in \cite{OUM10}.

\begin{lem}\label{lem:4.2} Let $(M,f,g)$ be a special matrix
  representation of a lagrangian $\bK_\sigma$-chain group $L$ on $V$. Let
  $f'$ be a $\bK_\sigma$-chain on $V$ such that $f'(x)\in \{c^*, c_*\mid c\in
  \bF^*\}$ for all $x\in V$. Then, $f'$ is eulerian if and only if
  $M[X]$ is non-singular with $X:=\{x\in V\mid f'(x)\ne c\cdot
  f(x)$ for some $c\in \bF^*\}$.
\end{lem}

\begin{pf*}{Proof.} (Proof already present in \cite{OUM10}.) Let
  $\{f_x\mid x\in V\}$ be the basis of $L$ associated with $f$ and $g$
  from Proposition \ref{prop:4.1}.  For each $y\in
  X$, there exists $d_y\in \bF^*$ such that
  \begin{align*}
    f'(y) &= \begin{cases} d_y\cdot f(y) & \textrm{if $y\notin
        X$},\\  d_y\cdot g(y) & \textrm{if $y\in
        X$}. \end{cases} \end{align*} 

  Assume that $M[X]$ is non-singular and let $h\in L$ such that
  $\mathbf{b}_\sigma(h(y),f'(y))=0$ for all $y\in V$. Let
  $h=\sum\limits_{z\in V} c_z\cdot f_z$. For all $y\notin X$, we have 
  \begin{align*}
\mathbf{b}_\sigma(h(y),f'(y)) &= 
\mathbf{b}_\sigma\left(\sum\limits_{z\in V}\left(c_z\cdot m_{zy}\cdot
f(y)\right)+c_y\cdot g(y), d_y\cdot f(y)\right) \\ &= \mathbf{b}_\sigma(c_y\cdot
g(y), d_y\cdot f(y)) \\ &= - c_y\cdot \sigma(d_y)\cdot
\epsilon(y)\cdot \sigma(1).
  \end{align*}
Hence, $c_y=0$ for all $y\notin X$. If $y\in X$, then
  \begin{align*}
\mathbf{b}_\sigma(h(y),f'(y)) & =
\mathbf{b}_\sigma\left(\sum\limits_{z\in X}\left(c_z\cdot m_{zy}\cdot
f(y)\right) + c_y\cdot g(y), d_y\cdot g(y)\right) \\ &=
\sum\limits_{z\in X} \left(c_z \cdot m_{z y}\cdot
\frac{\sigma(d_y)}{\sigma(1)}\cdot \mathbf{b}_\sigma(f(y),g(y))\right)
\\ & = (\sigma(d_y)\cdot \epsilon(y))\cdot \sum\limits_{z\in X} c_z
\cdot m_{z y}. \end{align*} For $\mathbf{b}_\sigma(h(y),f'(y))$ to
  being $0$, we must have $\sum\limits_{z\in X} \left(c_z \cdot m_{z
    y}\right) = 0$.  But, since $M[X]$ is non-singular, we have
  $\sum\limits_{z\in X} \left(c_z \cdot m_{z y}\right) = 0$ for all
  $y\in X$ if and only if $c_z=0$ for all $z\in X$. Therefore, we have
  $h=0$, \ie, $f'$ is eulerian.

  Assume now that $M[X]$ is singular. Hence, there exist $c_z$
  for $z\in X$, not all zero, such that for all $y\in X$, $\sum\limits_{z\in X}
  \left(c_z\cdot  m_{zy}\right) = 0$. Let $h:=\sum\limits_{z\in X}
  c_z\cdot 
  f_z$, which is not zero. Hence, for each $y\notin X$,
  \begin{align*}
    \mathbf{b}_\sigma(h(y),f'(y)) &= \frac{\sigma(d_y)}{\sigma(1)} \cdot
    \mathbf{b}_\sigma\left(\sum\limits_{z\in X} \left(c_z \cdot f_z (y)\right),
    f(y)\right) \\ & = \frac{\sigma(d_y)}{\sigma(1)} \cdot \left(
    \sum\limits_{z\in X} \left(c_z\cdot  m_{z y} \cdot
    \mathbf{b}_\sigma(f(y),f(y))\right) \right) = 0
  \end{align*}
  For each $y\in X$,  
  \begin{align*}
    \mathbf{b}_\sigma(h(y),f'(y)) & = \frac{\sigma(d_y)}{\sigma(1)}\cdot
    \mathbf{b}_\sigma\left(\sum\limits_{z\in X} \left(c_z \cdot f_z (y)\right),
    g(y)\right)\\ & = \frac{\sigma(d_y)}{\sigma(1)} \cdot \left(
    \sum\limits_{z\in X} \left(c_z\cdot  m_{zy} \cdot
    \mathbf{b}_\sigma(f(y),g(y))\right) \right) \\ &= \sigma(d_y)\cdot
    \epsilon(y) \cdot \left( \sum\limits_{z\in X} c_z\cdot  m_{z
      y} \right) = 0
  \end{align*}
  Since $h$ is not zero and $\mathbf{b}_\sigma(h(y),f'(y))=0$ for all
  $y\in V$, $f'$ is not eulerian. \qed
\end{pf*}

We now relate special matrix representations of a lagrangian
$\bK_\sigma$-chain group with the ones of its $\alpha\beta$-minors.

\begin{lem}\label{lem:4.1} Let $\{\alpha,\beta\}\subseteq
\{c^*,c_*\mid c\in \bF^*\}$ be minor-compatible.  Let $(M,f,g)$ be a
special matrix representation of a lagrangian $\bK_\sigma$-chain group
$L$ on $V$, and let $x\in V$. Then, $(M[V\setminus x],\restriction{f}{
  (V\setminus x)}, \restriction{g}{ (V\setminus x)})$ is a special
matrix representation of $L\minor{\alpha} x$ if $f(x)=c\cdot \alpha$,
otherwise of $L\minor{\beta}x$.
\end{lem}

\begin{pf*}{Proof.} We can assume by symmetry that
  $f(x)=c\cdot \alpha$. Let $\{f_x\mid x\in V\}$ be the basis of $L$
  associated with $f$ and $g$ from Proposition \ref{prop:4.1}.

  For all $y\in V\setminus x$, we have $f_y(x)=m_{yx}\cdot c\cdot
  \alpha$. Hence, $f_y\in L\minor{\alpha} x$ for all $y\in V\setminus
  x$. We claim that the set $\{\restriction{f_y}{(V\setminus x)}\mid
  y\in V\setminus x\}$ is linearly independent. Suppose the contrary
  and let $h:=\sum\limits_{y\in V\setminus x} c_y \cdot f_y \in L$ with
  $\restriction{h}{ (V\setminus x)} = 0$. Hence, $h(x) =
  \sum\limits_{y\in V\setminus x} \left(c_y\cdot m_{yx}\cdot c\cdot
  \alpha\right)$ and $h(y)=0$ for all $y\in V\setminus x$. Therefore,
  $\mathbf{b}_\sigma(h(z),f(z))= 0$ for all $z\in V$, contradicting
  the eulerian of $f$. By Proposition \ref{prop:3.4}, $L\minor{\alpha}
  x$ is lagrangian, \ie, $\dim(L\minor{\alpha} x)=|V\setminus x|$,
  hence $\{\restriction{f_y}{ (V\setminus x)}\mid y\in V\setminus
  x\}$ is a basis for $L\minor{\alpha} x$. But, this is actually the
  basis of $(M[V\setminus x],\restriction{f}{ (V\setminus x)},
  \restriction{g}{ (V\setminus x)})$ from Proposition
  \ref{prop:4.1}. \qed


\end{pf*}

We have then the following. 

\begin{prop}\label{prop:4.6} Let $\{\alpha,\beta\}\subseteq
\{c^*,c_*\mid c\in \bF^*\}$ be minor-compatible. Let $L$ and $L'$ be lagrangian
  $\bK_\sigma$-chain groups on $V$ and $V'$ respectively. Let
$(M,f,g)$ and $(M',f',g')$ be
  special matrix representations of $L$ and $L'$ respectively with
  $f(x):=\pm\alpha,\ g(x):=\beta$ for all $x\in V$, and
  $f'(x):=\pm\alpha,\ g'(x):=\beta$ for all $x\in V'$. If
  $L'=L\minor{\beta}X\minor{\alpha} Y$, then $M'= \big((M/M[A])[V']\big)\cdot
  I_Z$ with $A\subseteq X$ and $Z:=\{x\in V'\mid f'(x)=-f(x)\}$.
\end{prop}

\begin{pf*}{Proof.} If $X=\emptyset$, then by
  Lemma \ref{lem:4.1} $(M[V'],\restriction{f}{ V'}, \restriction{g}{ V'})$ is a special
  matrix representation of $L'$. By hypothesis, $g'=\restriction{g}{ V'}$. If
  we let $Z:=\{x\in V'\mid f'(x)=-f(x)\}$, then by Proposition
  \ref{prop:4.5} $(M[V']\cdot I_Z, f',g')$ is a special matrix
  representation of $L'$. Therefore, $M'=M[V']\cdot I_Z$ by
  Proposition \ref{prop:4.3}. We can now assume that $X\ne \emptyset$
  and is minimal with the property that there exists $Y$ such that
  $L'=L\minor{\beta}X\minor{\alpha} Y$.

  We claim that $M[X]$ is non-singular. Assume the contrary and let
  $f_1$ be the $\bK_\sigma$-chain on $V$ where $f_1(x)=f(x)$ if $x\notin X$,
  and $f_1(x)=g(x)$ otherwise. By Lemma \ref{lem:4.2}, $f_1$ is not
  eulerian. Hence, there exists $h\in L$ a non-zero $\bK_\sigma$-chain on $V$ such that
  $\mathbf{b}_\sigma(h(x),f_1(x)) = 0$ for all $x\in V$. Then, $\restriction{h}{
    V'} \in L'$. And since $\restriction{f_1}{ V'}=\restriction{f}{ V'} = f'$, we have
  $\restriction{h}{ V'} = 0$ ($f'$ is eulerian).  Moreover, there exists $z\in
  X$ such that $h(z) \ne 0$, otherwise it contradicts the fact that
  $f$ is eulerian (recall that for all $y\in V\setminus
  X,\ f_1(y)=f(y)$). By Lemma \ref{lem:3.2}, we have $h(z)=c_z\cdot
  \beta$, $c_z\in \bF^*$. Let $h'\in L$ such that $\restriction{h'}{ V'}\in
  L'$. Then, $\mathbf{b}_\sigma(h'(z), \beta) = 0$, and hence
  $\mathbf{b}_\sigma(h(z), h'(z)) = 0$. Thus by Lemma \ref{lem:3.2},
  $h'(z) = c_{h'}\cdot h(z)$.  Hence, $\restriction{(h'-c_{h'}\cdot h)}{ V'}
  \in L\minor{\beta}(X\setminus z)\minor{\alpha} (Y\cup z)$. But, we
  have $\restriction{(h'-c_{h'}\cdot h)}{ V'} = \restriction{h'}{ V'}$ because $\restriction{h}{
    V'} =0$. Therefore, $L\minor{\beta}X\minor{\alpha} Y \subseteq
  L\minor{\beta}(X\setminus z)\minor{\alpha} (Y\cup z)$. By
    Proposition \ref{prop:3.4}, $\dim(L\minor{\beta}X\minor{\alpha} Y)
    = |V'|$ and $\dim(L\minor{\beta}(X\setminus z)\minor{\alpha} (Y\cup
      z)) = |V\setminus (X\setminus z)\setminus (Y\cup
      z)|=|V'|$. Hence, $L\minor{\beta}X\minor{\alpha} Y =
      L\minor{\beta}(X\setminus z)\minor{\alpha} (Y\cup z)$. This
        contradicts the assumption that $X$ is minimal. Hence, $M[X]$
        is non-singular.

  Let $M_1:=P_X\cdot (M*X)$. By Proposition \ref{prop:4.4}, there
  exist $f_2$ and $g_2$ such that $L=(M_1,f_2,g_2)$. By Lemma
  \ref{lem:4.1}, $(M_1[V\setminus X], \restriction{f_2}{ V\setminus
    X}, \restriction{g_2}{ V\setminus X})$ is a matrix representation
  of $L\minor{\beta}X$. Notice that $\restriction{f_2}{ V\setminus X}
  = \restriction{f}{ V\setminus X}$ and $\restriction{g_2}{
    V\setminus X} = \restriction{g}{ V\setminus X}$. By Lemma
  \ref{lem:4.1}, $(M_1[V'], \restriction{f}{ V'}, \restriction{g}{
    V'})$ is a special matrix representation of
  $L\minor{\beta}X\minor{\alpha} Y$. But, $f'=\pm \restriction{f}{V'}$
  and $g'=\restriction{g}{ V'}$. Let $Z:=\{x\in V'\mid f'(x)=-
  f(x)\}$. By Proposition \ref{prop:4.5}, $(M_1[V']\cdot I_Z, f',g')$
  is a special matrix representation of $L'$. Therefore,
  $M'=M_1[V']\cdot I_Z$ by Proposition \ref{prop:4.3}. And, the fact
  that $M_1[V']=(M/M[X])[V']$ finishes the proof. \qed
\end{pf*}

We are now ready to prove the principal result of the paper. 

\begin{thm}\label{thm:4.2} Let $\bF$ be a finite field and $k$ a
  positive integer. For every infinite sequence $M_1,M_2,\ldots$ of
  $(\sigma_i,\epsilon_i)$-symmetric $(V_i,V_i)$-matrices over $\bF$ of
  $\bF$-rank-width at most $k$, there exist $i<j$ such that $M_i$ is
  isomorphic to $\big((M_j/M_j[A])[V']\big)\cdot I_Z$ with $A\subseteq
  V_j\setminus V'$ and $Z\subseteq V'$.
\end{thm}

\begin{pf*}{Proof.} Let $\alpha:= c^*$ and $\beta:=\widetilde{c^*}$
  for some $c\in \bF^*$. Since the set of sesqui-morphisms over $\bF$
  is finite, we can assume by taking a sub-sequence that each matrix
  $M_i$ is $(\sigma,\epsilon_i)$-symmetric, for some sesqui-morphism
  $\sigma:\bF\to \bF$. For each $i$, let $f_i$ and $g_i$ be
  $\bK_\sigma$-chains on $V_i$ with $f_i(x):=\epsilon_i(x)\cdot
  \alpha$ and $g_i(x):=\beta$ for all $x\in V_i$. Let $L_i$ be
  $(M_i,f_i,g_i)$. By Theorem \ref{thm:3.3}, there exist $i<j$ such
  that $L_i$ is simply isomorphic to an $\alpha\beta$-minor of
  $L_j$. Let $X,Y\subseteq V_j$ such that $L_i$ is simply isomorphic
  to $L_j\minor{\beta}X\minor{\alpha} Y$. Let $V':=V_j\setminus (X\cup
  Y)$. By Proposition \ref{prop:4.6}, $M_i$ is isomorphic to
  $\big((M_j/M_j[A])[V']\big)\cdot I_Z$ with $A\subseteq X$ and $Z\subseteq
  V'$.\qed
\end{pf*}

Since each symmetric (or skew-symmetric) $(V,V)$-matrix is a
$(\sigma,\epsilon)$-symmetric $(V,V)$-matrix with $\epsilon(x)=1$ for
all $x\in V$, and $\sigma$ being symmetric (or skew-symmetric),
Theorem \ref{thm:2.1} is a corollary of Theorem \ref{thm:4.2}.  It is
worth noticing as noted in \cite{OUM10} that the well-quasi-ordering
results in \cite{GGW02,OUM05,RS90} are corollaries of Theorem
\ref{thm:2.1}, hence of Theorem \ref{thm:4.2}. We give some other
corollaries about graphs in the next section.

\section{Applications to Graphs}\label{sec:5}

Clique-width was defined by Courcelle et al. \cite{CER93} for graphs
(directed or not, with edge-colours or not). But, the notion of
rank-width introduced by Oum and Seymour in \cite{OS06} and studied by
Oum (see for instance \cite{OUM05,OUM08}) concerned only undirected
graphs. Rao and myself we generalised in \cite{KR11} the notion of
rank-width to directed graphs, and more generally to edge-coloured
graphs.  We give well-quasi-ordering theorems for directed graphs and
edge-coloured graphs.

\subsection{The Case of Edge-Coloured Graphs}

Let $C$ be a (possibly infinite) set that we call the \emph{colours}.
A \emph{$C$-coloured graph} $G$ is a tuple $(V_G, E_G, \ell_G)$ where
$(V_G,E_G)$ is a directed graph and $\ell_G: E_G \to
2^C\setminus\{\emptyset\}$ is a function. Its associated
\emph{underlying graph} $\supp{G}$ is the directed graph
$(V_G,E_G)$. Two $C$-coloured graphs $G$ and $H$ are isomorphic if
there is an isomorphism $h$ between $\supp{G}$ and $\supp{H}$ such
that for every $(x,y)\in E_G$, $\ell_G((x,y))=\ell_H((h(x),h(y))$.  We
call $h$ an \emph{isomorphism} between $G$ and $H$. It is worth
noticing that an edge-uncoloured graph can be seen as an edge-coloured
graph where all the edges have the same colour. 

The notion of rank-width of $C$-coloured graphs is based on the
$\bF$-rank-width of $(\sigma,\epsilon)$-symmetric matrices. Let $\bF$
be a field. An \emph{$\bF^*$-graph} $G$ is an $\bF^*$-coloured graph
where for every edge $(x,y)\in E_G$, we have $\ell_G((x,y))\in \bF^*$,
\ie, each edge has exactly one colour in $\bF^*$. It is clear that
every directed graph is an $\bF_2^*$-graph. One interesting point is
that every $\bF^*$-graph $G$ can be represented by a
$(V_G,V_G)$-matrix $\matg$ over $\bF$, that generalises the adjacency
matrix of directed graphs, such that
\begin{align*}
\matgind{x}{y}:= \begin{cases} \ell_G((x,y)) & \textrm{if $(x,y)\in
    E_G$},\\ 0 & \textrm{otherwise}.
\end{cases}
\end{align*}

If $\matg$ is $(\sigma,\epsilon)$-symmetric, we call $G$ a
\emph{$(\sigma,\epsilon)$-symmetric $\bF^*$-graph}. It is worth noticing that
in this case $\supp{G}$ is undirected. Not all $\bF^*$-graphs are
$(\sigma,\epsilon)$-symmetric, however we have the following.

\begin{prop}[\cite{KR11}]\label{prop:5.1} Let $\bF$ be a finite
  field. Then, one can construct a sesqui-morphism $\sigma:\bF^2\to
  \bF^2$ where $\bF^2$ is an algebraic extension of $\bF$ of order
  $2$. Moreover, for every $\bF^*$-graph $G$, one can associate a
  $\sigma$-symmetric $(\bF^2)^*$-graph
  $\widetilde{G}$ such that for every $\bF^*$-graphs $G$ and $H$,
  $\widetilde{G}$ and $\widetilde{H}$ are isomorphic if and only if
  $G$ and $H$ are isomorphic. 
\end{prop}

In order to define a notion of rank-width for $C$-coloured graphs, we
proceed as follows. For a $C$-coloured graph $G$, let $\Pi(G)\subseteq
2^C$ be the set of subsets of $C$ appearing as colours of edges in
$G$.

\begin{enumerate}
\item take an injection $i:\Pi(G) \to \bF^*$ for a large enough
  finite field $\bF$ and let $G'$ be the $\bF^*$-graph
  obtained from $G$ by replacing each edge colour $A\subseteq C$ by
  $i(A)$. If the $\bF^*$-graph $G'$ is $(\sigma,\epsilon)$-symmetric
  for some sesqui-morphism $\sigma:\bF\to \bF$, then define the
  $\bF$-rank-width of $G$ as the $\bF$-rank-width of $\mat{G'}$. Otherwise, 
\item take $\widetilde{G'}$ from
  Proposition \ref{prop:5.1}. $\mat{\widetilde{G'}}$ is
  $\sigma$-symmetric for some $\sigma:\bF^2\to \bF^2$. The
  \emph{$\bF^2$-rank-width} of $G$ will be defined as the
  $\bF^2$-rank-width of $\mat{\widetilde{G'}}$.
\end{enumerate}

The choice of the injection in step (1) above is not unique and leads
to different representations of $C$-coloured graphs, and then
different parameters. However, as proved in \cite{KR11}, the
parameters are equivalent. Therefore, in order to investigate the
structure of $C$-coloured graphs, we can concentrate our efforts in
$(\sigma,\epsilon)$-symmetric $\bF^*$-graphs. The authors in
\cite{KR11} did only consider  $\sigma$-symmetric graphs. We relax
this constraint because we may have some $\bF^*$-graphs which are
$(\sigma,\epsilon)$-symmetric but are not $\sigma'$-symmetric at all,
for all sesqui-morphisms $\sigma':\bF\to \bF$. Examples of such graphs
are $\bF^*$-graphs $G$ where $\matg$ is obtained from a
$\sigma$-symmetric matrix by multiplying some rows and/or columns by
$-1$. 

All the results, but the well-quasi-ordering theorem, concerning the
rank-width of undirected graphs are generalised in \cite{KR11} to the
$\bF$-rank-width of $\sigma$-symmetric loop-free $\bF^*$-graphs. These
results extend easily to $(\sigma,\epsilon)$-symmetric
$\bF^*$-graphs. We prove here two well-quasi-ordering theorems for
$(\sigma,\epsilon)$-symmetric $\bF^*$-graphs. For that, we will
derive from the principal pivot transform two notions of pivot-minor:
one that preserves the loop-freeness and one that does not. 

We recall that a pair $(P,Q)$ of non-singular diagonal $(V,V)$-matrices is
$\sigma$-compatible if $p_{xx}^{-1}=\sigma(q_{xx})\cdot
\sigma(1)^{-1}$ (equivalently $q_{xx}^{-1}=\sigma(p_{xx})\cdot
\sigma(1)^{-1}$) for all $x\in V$, and for $X \subseteq V$, 
 $P_X$ and $I_X$ are the non-singular diagonal $(V,V)$-matrices
where
\begin{align*}
 P_X[x,x] &:= \begin{cases} \sigma(-1) & \textrm{if $x\in X$},\\ 1 &
  \textrm{otherwise}, \end{cases} & \quad \textrm{and} & \quad I_X[x,x]
 := \begin{cases} -1 & \textrm{if $x\in X$},\\ 1 &
  \textrm{otherwise.} \end{cases}
\end{align*}

\begin{defn}[$\sigma$-loop-pivot complementation]\label{defn:5.2} Let $G$ be a
  $(\sigma,\epsilon)$-symmetric $\bF^*$-graph and let $X\subseteq V_G$
  such that $\matg[X]$ is non-singular. An $\bF^*$-graph $G'$ is a \emph{$\sigma$-loop-pivot
    complementation of $G$ at $X$} if $\mat{G'}:=I_Z\cdot P\cdot P_X\cdot (M*X)
  \cdot Q^{-1} \cdot I_{Z'}$ for some  $Z,Z'\subseteq V_G$, and $(P,Q)$ 
  a pair of $\sigma$-compatible diagonal $(V_G,V_G)$-matrices. 

  An $\bF^*$-graph $G'$ is \emph{$\sigma$-loop-pivot equivalent} to
  $G$ if $G'$ is obtained from $G$ by applying a sequence of $\sigma$-loop-pivot
    complementations.  An $\bF^*$-graph $H$ is a
  $\sigma$-loop-pivot-minor of  $G$ if $H$ is isomorphic to $G'[V'],\ V'\subseteq
  V_G$, where $G'$ is $\sigma$-loop-pivot equivalent to $G$.
\end{defn} 

The $\sigma$-loop-pivot complementation does not clearly preserve the
loop-freeness. A corollary of Theorem \ref{thm:4.1},
and Propositions \ref{prop:4.4}, \ref{prop:4.5} and \ref{prop:4.7} is
the following.

\begin{cor}\label{cor:5.1} \begin{enumerate} 
\item Let $G$ be a $(\sigma,\epsilon)$-symmetric $\bF^*$-graph. If
  $G'$ is $\sigma$-loop-pivot equivalent to $G$, then $G'$ is
  $(\sigma,\epsilon')$-symmetric for some $\epsilon':V_G\to
  \{+1,-1\}$.

\item Let $G$ and $G'$ be respectively
  $(\sigma,\epsilon)$ and $(\sigma,\epsilon')$-symmetric
  $\bF^*$-graphs. If $G'$ is $\sigma$-loop-pivot equivalent to $G$, then
  $\frwd{G'}=\frwd{G}$. If $G'$ is a $\sigma$-loop-pivot-minor of $G$, then
  $\frwd{G'}\leq \frwd{G}$.
\end{enumerate}
\end{cor}

We now introduce a variant of the $\sigma$-loop-pivot complementation
that preserves the loop-freeness and prove that Corollary
\ref{cor:5.1} still holds.

\begin{defn}[$\sigma$-pivot complementation]\label{defn:5.1} Let $G$ be a
  $(\sigma,\epsilon)$-symmetric loop-free $\bF^*$-graph and let
  $X\subseteq V_G$ such that $\matg[X]$ is non-singular. A loop-free
  $\bF^*$-graph $H$ is a \emph{$\sigma$-pivot complementation of $G$
    at $X$} if $\mat{H}$ is obtained from $\mat{G'}$, $G'$ a
  $\sigma$-loop-pivot complementation of $G$ at $X$, by replacing each
  diagonal entry by $0$.

  A loop-free $\bF^*$-graph $G'$ is \emph{$\sigma$-pivot equivalent} to
  $G$ if $G'$ is obtained from $G$ by applying a sequence of $\sigma$-pivot
    complementations.  A loop-free $\bF^*$-graph $H$ is a
  $\sigma$-pivot-minor of  $G$ if $H$ is isomorphic to $G'[V'],\ V'\subseteq
  V_G$, where $G'$ is $\sigma$-pivot equivalent to $G$.
\end{defn} 

It is clear that the $\sigma$-pivot complementation preserves the
loop-freeness. The proof of the following is straightforward.

\begin{prop}\label{prop:5.2} Let $(M,f,g)$ be a matrix representation
  of a lagrangian $\bK_\sigma$-chain group $L$ on $V$ and let $M'$ be obtained
  from $M$ by replacing each diagonal entry by $0$. Let $g'$
  be the $\bK_\sigma$-chain on $V$ with $g'(x):=m_{xx}\cdot f(x) +
  g(x)$. Then, $(M',f,g')$ is a matrix representation of $L$.  
\end{prop}

The following is hence true. 

\begin{cor}\label{cor:5.2} \begin{enumerate} 
\item Let $G$ be a $(\sigma,\epsilon)$-symmetric loop-free $\bF^*$-graph. If
  $G'$ is $\sigma$-pivot equivalent to $G$, then $G'$ is
  $(\sigma,\epsilon')$-symmetric for some $\epsilon':V_G\to
  \{+1,-1\}$.

\item Let $G$ and $G'$ be respectively $(\sigma,\epsilon)$ and
  $(\sigma,\epsilon')$-symmetric loop-free $\bF^*$-graphs. If $G'$ is
  $\sigma$-pivot equivalent to $G$, then $\frwd{G'}=\frwd{G}$. If $G'$
  is a $\sigma$-pivot-minor of $G$, then $\frwd{G'}\leq \frwd{G}$.
\end{enumerate}
\end{cor}

As corollaries of Theorem \ref{thm:4.2}, we have the following
well-quasi-ordering theorems for $\bF^*$-graphs. 

\begin{thm}\label{thm:5.1} Let $\bF$ be a finite field and $k$ a
  positive integer. For every infinite sequence $G_1,G_2,\ldots$ of
  $(\sigma_i,\epsilon_i)$-symmetric $\bF^*$-graphs of
  $\bF$-rank-width at most $k$, there exist $i<j$ such that $G_i$ is
  isomorphic 
 a $\sigma$-loop-pivot-minor of $G_j$.
\end{thm}

\begin{pf*}{Proof.} Let $\mat{G_1},\mat{G_2}, \ldots$ be the infinite
  sequence of $(\sigma_i,\epsilon_i)$-symmetric
  $(V_{G_i},V_{G_i})$-matrices over $\bF$ associated with the infinite
  sequence $G_1,G_2,\ldots$. By definition,
  $\frwd{G_i}=\frwd{\mat{G_i}}$. From Theorem \ref{thm:4.2}, there
  exist $i<j$ such that $\mat{G_i}$ is isomorphic to
  $\big((\mat{G_j}/\mat{G_j}[A])[V']\big)\cdot I_Z$ with $A,V',Z\subseteq
  V_{G_j}$. But, that means that $G_i$ is isomorphic to a
  $\sigma$-loop-pivot-minor of $G_j$. \qed
\end{pf*}

\begin{thm}\label{thm:5.2} Let $\bF$ be a finite field and $k$ a
  positive integer. For every infinite sequence $G_1,G_2,\ldots$ of
  $(\sigma_i,\epsilon_i)$-symmetric loop-free $\bF^*$-graphs of
  $\bF$-rank-width at most $k$, there exist $i<j$ such that $G_i$ is
  isomorphic to 
 a $\sigma$-pivot-minor of $G_j$.
\end{thm}

\begin{pf*}{Proof.} Let
  $\mat{G_1},\mat{G_2}, \ldots$ be the infinite sequence of
  $(\sigma_i,\epsilon_i)$-symmetric $(V_{G_i},V_{G_i})$-matrices over
  $\bF$ associated with the infinite sequence $G_1,G_2,\ldots$. By
  definition, $\frwd{G_i}=\frwd{\mat{G_i}}$.  From Theorem
  \ref{thm:4.2}, there exist $i<j$ such that $\mat{G_i}$ is isomorphic
  to $((\mat{G_j}/\mat{G_j}[A])[V'])\cdot I_Z$ with $A,V',Z\subseteq
  V_{G_j}$. Since, $G_i$ is loop-free, this means that the diagonal
  entries of $\big((\mat{G_j}/\mat{G_j}[A])[V']\big)\cdot I_Z$ are equal to
  $0$. Hence, $(\mat{G_j}*A)[V']$ has only zero in its diagonal
  entries. Then, $G_i$ is isomorphic to a $\sigma$-pivot-minor of
  $G_j$. \qed
\end{pf*}

\subsection{A Specialisation to  Directed Graphs}

We discuss in this section a corollary about directed graphs. Let us
first recall the rank-width notion of directed graphs. We recall that
$\field{4}$ is the finite field of order four. We let
$\{0,1,\gfa,\gfb\}$ be its elements with the property that
$1+\gfa+\gfb=0$ and $\gfa^3=1$. Moreover, it is of characteristic
$2$. We let $\sigma_4:\field{4}\to \field{4}$ be the automorphism
where $\sigma_4(\gfa) = \gfb$ and $\sigma_4(\gfb) = \gfa$.  It is
clearly a sesqui-morphism.

For every directed graph $G$, let $\widetilde{G}:=(V_G,E_G\cup\{(y,x)\vert (x,y)\in E_G\},\ell_{\tG})$ be the
$\gfq^*$-graph where for every pair of vertices $(x,y)$:
\begin{align*}
  \ell_{\tG}((x,y)) &:= \begin{cases} 1 & \textrm{if $(x,y)\in
      E_G\ \textrm{and}\ (y,x)\in E_G$},\\ \gfa &
    \textrm{$(x,y)\in E_G\ \textrm{and}\ (y,x)\notin E_G$},\\ \gfb
    & \textrm{$(y,x)\in E_G\ \textrm{and}\ (x,y)\notin E_G$},\\ 0
    & \textrm{otherwise}.
  \end{cases}
\end{align*}

It is straightforward to verify that $\widetilde{G}$ is
$\sigma_4$-symmetric and there is a one-to-one correspondence between
directed graphs and $\sigma_4$-symmetric $\bF_4^*$-graphs.  The
\emph{rank-width} of a directed graph $G$, denoted by $\Qrwd{G}$, is
the $\gfq$-rank-width of $\widetilde{G}$ \cite{KR11}. One easily
verifies that if $G$ is an undirected graph, then the rank-width of
$G$ is exactly the $\bF_4$-rank-width of $\widetilde{G}$.

\begin{figure}[h!]
  \centering
    \input{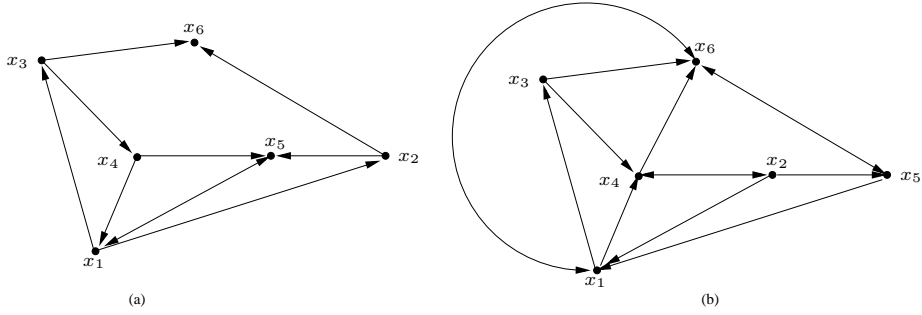}
  \caption{(a) A directed graph $G$. (b) The directed graph obtained
    after a pivot-complementation of $G$ at $\{x_2,x_5\}$. If you
    apply a loop-pivot-complementation of $G$ at $\{x_2,x_5\}$, you
    obtain the graph in (b) with a loop at $x_1$.}
  \label{fig:5.1}
\end{figure}

A directed graph $H$ is \emph{loop-pivot equivalent}
(resp. \emph{pivot equivalent}) to a directed graph $G$ if
$\widetilde{H}$ is $\sigma_4$-loop-pivot equivalent
(resp. $\sigma_4$-pivot equivalent) to $\widetilde{G}$; and $H$ is a
\emph{loop-pivot-minor} (resp. \emph{pivot-minor}) of $G$ if
$\widetilde{H}$ is a $\sigma_4$-loop-pivot minor
(resp. $\sigma_4$-pivot minor) of $\widetilde{G}$. Since there is a
one-to-one correspondence between $\sigma_4$-symmetric
$\bF_4^*$-graphs and directed graphs, loop-pivot equivalence
(resp. pivot-equivalence) and loop-pivot minor (resp. pivot-minor) are
well-defined in directed graphs. Figure \ref{fig:5.1} shows an
example of loop-pivot complementation and pivot complementation.

As a consequence of Theorems \ref{thm:5.1} and \ref{thm:5.2} we have
the following which generalises \cite[Theorem 4.1]{OUM08}. 

\begin{thm} \label{thm:5.3} Let $k$ be a positive
  integer. 
  \begin{enumerate} 
  \item For every infinite sequence $G_1,G_2,
    \ldots$ of directed graphs of rank-width at most
    $k$, there exist $i<j$ such that $G_i$ is isomorphic to a loop-pivot-minor of
    $G_j$.

  \item For every infinite sequence $G_1,G_2, \ldots$ of loop-free
    directed graphs of rank-width at most $k$, there exist $i<j$ such
    that $G_i$ is isomorphic to a pivot-minor of $G_j$.
  \end{enumerate}
\end{thm}

\section{Delta-Matroids and Chain Groups} \label{sec:6}

In this section we discuss some consequences of results in Sections
\ref{sec:3} and \ref{sec:4} about \emph{delta-matroids}. If $V$ is a
finite set, then $\cF\subseteq 2^V$ is said to satisfy the
\emph{symmetric exchange axiom} if:
\begin{quote}
  (SEA)\ for $F,F'\in \cF$, for $x\in F\triangle F'$, there exists
  $y\in F'\triangle F$ such that $F\triangle \{x,y\}\in \cF$. 
\end{quote}

A \emph{set system} is a pair $(V,\cF)$ where $V$ is finite and
$\emptyset \ne \cF \subseteq 2^V$.  A \emph{delta-matroid} is a
set-system $(V,\cF)$ such that $\cF$ satisfies (SEA); the elements of
$\cF$ are called \emph{feasible sets}. Delta-matroids were
introduced in \cite{BOU87b}, and as for matroids, are characterised
by the validity of a greedy algorithm. We
recall that a set system $\cM:=(V,\cB)$ is a \emph{matroid} if $\cB$,
called the set of \emph{bases}, satisfy the following \emph{Exchange
  Axiom}
\begin{quote}
  (EA) for $B,B'\in \cB$, for $x\in B\setminus B'$, there exists
  $y\in B'\setminus B$ such that $B\triangle \{x,y\}\in \cB$.
\end{quote}

It is worth noticing that a matroid is also a delta-matroid (see
\cite{BOU87b,BOU88,GEE95} for other examples of delta-matroids). 

For a set system $\cS=(V,\cF)$ and $X\subseteq V$, we let
$\cS\triangle X$ be the set system $(V,\cF\triangle X)$ where
$\cF\triangle X:=\{F\triangle X\mid F\in \cF\}$. We have that
$\cF\triangle X$ satisfies (SEA) if and only if $\cF$ satisfies
(SEA). Hence, $\cS$ is a delta-matroid if and only if $\cS\triangle X$
is. A delta-matroid $\cS=(V,\cF)$ is said \emph{equivalent} to a
delta-matroid $\cS'=(V,\cF')$ if there exists $X\subseteq V$ such that
$\cS=\cS'\triangle X$.  If $M$ is a $(V,V)$-matrix over a field $\bF$,
we let $\cS(M)$ be the set system $(V,\cF(M))$ where
$\cF(M):=\{X\subseteq V\mid M[X]$ is non-singular$\}$. The following
is due to Bouchet \cite{BOU88}.

\begin{thm}[\cite{BOU88}]\label{thm:6.1} Let $M$ be a 
  matrix over $\bF$ of symmetric type, \ie, $M$ is
  $(\sigma,\epsilon)$-symmetric with $\sigma$ (skew) symmetric. Then, $\cS(M)$ is
  a delta matroid.
\end{thm}

Delta-matroids equivalent to $\cS(M)$, for some matrix $M$ over
$\bF$ of symmetric type, are called \emph{representable over $\bF$}
\cite{BOU88}.  A slight modification of the proof given in
\cite{GEE95} extends Theorem \ref{thm:6.1} to all
$(\sigma,\epsilon)$-symmetric matrices.

\begin{thm}\label{thm:6.2} Let $M$ be a $(\sigma,\epsilon)$-symmetric
  $(V,V)$-matrix over $\bF$. Then, $\cS(M)$ is a delta matroid.
\end{thm}

Let us recall the following from Tucker.

\begin{thm}[\cite{TUC60}] \label{thm:6.3} Let $M$ be a $(V,V)$-matrix such that
  $M[X]$ is non-singular. For any $Z\subseteq V$, we have \begin{align*}
    \det((M{*}X)[Z]) & = \pm
    \frac{\det(M[Z\triangle
        X])}{\det(A)}. \end{align*}
\end{thm}

\begin{pf*}{Proof of Theorem \ref{thm:6.2}.} Let $X,Y\subseteq V$ such
  that $M[X]$ and $M[Y]$ are non-singular. Let $x\in X\triangle Y$.
  Let $M':=P_X\cdot (M*X)$.  By Theorem \ref{thm:6.3}, $M'[Z]$ is
  non-singular if and only if $M[Z\triangle X]$ is non-singular.
  Assume $m'_{xx}\ne 0$, then if we take $y:=x$, we have that
  $M[X\triangle \{x\}]$ is non-singular. Suppose that
  $m'_{xx}=0$. Since $M'[X\triangle Y]$ is non-singular, there exists
  $y\in X\triangle Y$ such that $m'_{xy}\ne 0$ and because $M'$ is
  $(\sigma,\epsilon)$-symmetric, $m'_{yx} \ne 0$. Hence, $M'[\{x,y\}]$
  is non-singular, \ie, $M'[X\triangle \{x,y\}]$ is non-singular. \qed
\end{pf*}

A consequence of Theorem \ref{thm:6.2} is that we can extend the
notion of representability of delta-matroids by the following.
\begin{quote}
  A delta-matroid is \emph{representable over $\bF$} if it is equivalent to
  $\cS(M)$ for some $(\sigma,\epsilon)$-symmetric matrix $M$ over $\bF$.
\end{quote}

It is worth noticing from Proposition \ref{prop:2.2} that over prime
fields this notion of representability is the same as the one defined
by Bouchet \cite{BOU88}. We now discuss some other corollaries. First,
if $M$ is a $(\sigma,\epsilon)$-symmetric $(V,V)$-matrix, then for any
$X\subseteq V$ such that $M[X]$ is non-singular, $\cS(M)\triangle X =
\cS(M')$ for any $M':=I_Z\cdot P\cdot P_X\cdot (M*X) \cdot Q^{-1}
\cdot I_{Z'}$ for some $Z,Z'\subseteq V$, and $(P,Q)$ a pair of
$\sigma$-compatible diagonal $(V,V)$-matrices.




Lemma \ref{lem:4.2} characterises non-singular principal submatrices
of $(\sigma,\epsilon)$-symmetric matrices in terms of eulerian
$\bK_\sigma$-chains of their associated lagrangian $\bK_\sigma$-chain
groups. One can derive from this a characterisation of representable
delta-matroids in terms of lagrangian $\bK_\sigma$-chain groups.  

One can derive from Theorem \ref{thm:4.2} a well-quasi-ordering
theorem for representable delta-matroids as follows. Let the
\emph{branch-width} of a delta-matroid $\cS$ representable over $\bF$
as $\min \{\frwd{M}\mid \cS(M)$ is equivalent to $\cS\}$. A
delta-matroid $\cS'$ is a \emph{minor} of a delta-matroid
$\cS=(V,\cF)$ if there exist $X,Y\subseteq V$ such that
$\cS'=(V\setminus (X\cup Y),\{(F\triangle X)\setminus Y\mid F\in
\cF\})$.  An extension of \cite[Theorem 7.3]{OUM10} is the following.

\begin{thm}\label{thm:6.4} Let $\bF$ be a finite field and $k$ a
  positive integer. Every infinite sequence $\cS_1,\cS_2,\ldots$ of
  delta-matroids representable over $\bF$ of branch-width at most $k$
  has a pair $i<j$ such that $\cS_i$ is isomorphic to a minor of
  $\cS_j$.
\end{thm}

\begin{pf*}{Proof.} Let $M_1,M_2,\ldots$ be
  $(\sigma_i,\epsilon_i)$-symmetric matrices  over $\bF$ such that, for every $i$, $\cS_i$ is equivalent
  to $\cS(M_i)$ and the branch-width of $\cS_i$ is equal to the
  $\bF$-rank-width of $M_i$. By Theorem \ref{thm:4.2}, there exist $i<j$
  such that $M_i$ is isomorphic to $(M_j/M_j[A])[V']\cdot I_{Z}$ with
  $A\subseteq V_j\setminus V'$ and $Z\subseteq V'\subseteq V_j$. Hence, $\cS_i$ is
  isomorphic to a minor of $\cS_j$. \qed
\end{pf*}

We conclude by some questions.  It is well-known that columns of a
matrix over a field yields a matroid. It would be challenging to
characterise matrices whose non-singular principal submatrices yield a
delta-matroid.  Currently, there is no connectivity function for
delta-matroids. Another challenge is to find a connectivity function
for delta-matroids that subsumes the connectivity function of matroids
and such that if a delta-matroid is equivalent to $\cS(M)$, then the
branch-width of $\cS(M)$ is proportional to the $\bF$-rank-width of
$M$.

\begin{ack} We would like to thank S. Oum for letting at our
  disposal a first draft of \cite{OUM10}, which was of great help for
  our understanding of the problem. We thank also B. Courcelle and the
  anonymous referee for their helpful comments. The author is
  supported by the DORSO project of ``Agence Nationale Pour la
  Recherche''.
\end{ack}

  \end{document}